%% file: main.tex
\documentclass[11pt]{article}
\usepackage[a4paper, margin=3cm]{geometry}

\usepackage{tikz}
\usepackage{enumitem}
\usetikzlibrary{decorations.markings,arrows.meta}
\usepackage{xurl}
\usepackage{amsthm}
\usepackage{amsmath}
\newcommand{\E}{\mathbb{E}}
\renewcommand{\P}{\mathbb{P}}
\DeclareMathOperator{\Cov}{Cov}
\DeclareMathOperator{\Var}{Var}
\usepackage{hyperref}
\usepackage[nameinlink]{cleveref}
\usepackage{amsfonts}
\usepackage[style=numeric,maxbibnames=99]{biblatex}
\usepackage{epsfig}
\usepackage{float}
\usepackage{stmaryrd}
\usepackage{graphicx}
\graphicspath{ {./drawings/} }

\usepackage{subfiles}

\newtheorem{theorem}{Theorem}[section]
\newtheorem{conjecture}[theorem]{Conjecture}
\crefname{conjecture}{conjecture}{conjectures}

\newtheorem{lemma}[theorem]{Lemma}
\theoremstyle{definition}
\newtheorem{definition}[theorem]{Definition}

\newcommand{\unlab}[1]{\left\llbracket #1\right\rrbracket}

\title{Seeking the dense half in $K_4$-free graphs \thanks{This work was supported by the National Science Centre grant 2021/42/E/ST1/00193.}}
\author{Ignacy Buczek\thanks{Faculty of Mathematics and Computer Science, Jagiellonian University, {\L}ojasiewicza 6, 30-348 Krak\'{o}w, Poland. E-mail: {\tt ignacy.buczek@doctoral.uj.edu.pl}.} \and
Andrzej Grzesik\thanks{Faculty of Mathematics and Computer Science, Jagiellonian University, {\L}ojasiewicza 6, 30-348 Krak\'{o}w, Poland. E-mail: {\tt andrzej.grzesik@uj.edu.pl}.} \and
Piotr Kuc\thanks{Faculty of Mathematics and Computer Science, Jagiellonian University, {\L}ojasiewicza 6, 30-348 Krak\'{o}w, Poland. E-mail: {\tt piotr.kuc@student.uj.edu.pl}.}}
\date{\today}

\begin{document}

\maketitle

\begin{abstract}
We show that every $K_4$-free graph on $n$ vertices has a balanced bipartition such that neither part contains more than $0.0671n^2 + o(n^2)$ edges.
This improves the previous best bound of $0.074n^2$ by Xu and Zhang.
We also formulate a conjecture that the tight bound for this problem is $\frac{425}{6498}n^2 \approx 0.0654n^2$, and present a construction showing that it would be best possible.
\end{abstract}

\input{utils/flag-utils.tex}
\input{flags.tex}

\section{Introduction}
\input{sections/introduction.tex}

\section{Lower bound}\label{sec:lower-bound}
\input{sections/lower-bound.tex}

\section{Flag algebras}\label{sec:flag-algebras}
\input{sections/flag-algebras.tex}

\section{Upper bound}\label{sec:upper-bound}
\input{sections/upper-bound.tex}

\section{Acknowledgements}\label{sec:acknowledgement}
\input{sections/acknowledgement.tex}

\section{AI statement}\label{sec:ai}
\input{sections/ai.tex}

\printbibliography
\end{document}

%% file: utils/flag-utils.tex
\tikzset{flag/.style={scale=0.30}}

\tikzset{default/.style={
    draw,
    minimum size=1.8mm,
    inner sep=0pt,
    font=\tiny,
  }
}

\tikzset{root/.style={
    rectangle,
  }
}

\tikzset{vertex/.style={
    circle,
    fill=black,
  }
}

\tikzset{red/.style={
    fill=red,
    text=white
  }
}

\tikzset{blue/.style={
    fill=blue,
    text=white
  }
}

\tikzset{cyan/.style={
    fill=cyan,
    text=white
  }
}

\tikzset{white/.style={
    fill=white,
    text=black
  }
}

\tikzset{black/.style={
    fill=black,
    text=white
  }
}

\def\FlagInit{
  \pgfdeclarelayer{roots}
  \pgfdeclarelayer{vertices}
  \pgfdeclarelayer{edges}
  \pgfsetlayers{edges, vertices, roots}
}

\def\VertexCoordinates#1#2{
  \pgfmathsetmacro{\N}{#1 - 1}
  \foreach \i in {0, ..., \N} {
    \pgfmathsetmacro{\angle}{360/#1 * \i + #2}%
    \pgfmathsetmacro{\x}{cos(\angle)}%
    \pgfmathsetmacro{\y}{sin(\angle)}%
    \coordinate (x\i) at (\x,\y);%
  }
}

\def\RootVertex#1#2#3{
  \begin{pgfonlayer}{roots}
    \node[default, root, #3] at (x#1) {\the\numexpr#2+1\relax};
  \end{pgfonlayer}
}

\def\DrawVertex#1#2{
  \begin{pgfonlayer}{vertices}
    \node[default, vertex, #2] at (x#1) {};
  \end{pgfonlayer}
}

\def\RootVertices#1#2{
  \pgfmathsetmacro{\N}{#1 - 1}
  \foreach \i in {0, ..., \N} {
    \RootVertex{\i}{\i}{#2}
  }
}

\def\Connect#1#2#3{
  \begin{pgfonlayer}{edges}
    \draw[#3] (x#1) -- (x#2);
  \end{pgfonlayer}
}

\def\ConnectAny#1#2#3{
  \begin{pgfonlayer}{edges}
    \draw[#3] (x#1) -- (x#2)
      node[midway, circle, fill=black, inner sep=0.7pt] {};
  \end{pgfonlayer}
}

\def\ConnectOriented#1#2#3{
  \begin{pgfonlayer}{edges}
    \draw[#3, postaction={
      decorate,
      decoration={
        markings,
        mark=at position 0.73 with {\arrow[scale=0.9]{Latex}}
      }
    }] (x#1) -- (x#2);
  \end{pgfonlayer}
}

\def\Cycle#1{
  \pgfmathsetmacro{\N}{#1 - 1}
  \foreach \i in {0, ..., \N} {
    \pgfmathtruncatemacro{\nextvertex}{mod(\i + 1, #1)}
    \Connect{\i}{\nextvertex}{black}
  }
}

\newcommand{\vc}[1]{\ensuremath{\vcenter{\hbox{#1}}}}

%% file: flags.tex
\newcommand{\Vertex}[1]{
  \vc{
    \begin{tikzpicture}[flag]
      \FlagInit
      \VertexCoordinates{1}{0}
      \DrawVertex{0}{#1}
    \end{tikzpicture}
  }
}

\newcommand{\RootedVertex}[1]{
  \vc{
    \begin{tikzpicture}[flag]
      \FlagInit
      \VertexCoordinates{1}{0}
      \RootVertex{0}{0}{#1}
    \end{tikzpicture}
  }
}

\newcommand{\VertexM}[2]{
  \vc{
    \begin{tikzpicture}[flag]
      \FlagInit
      \VertexCoordinates{2}{240}
      \RootVertex{0}{0}{#1}
      \DrawVertex{1}{#2}
      \Connect{0}{1}{black}
    \end{tikzpicture}
  }
}

\newcommand{\VertexN}[2]{
  \vc{
    \begin{tikzpicture}[flag]
      \FlagInit
      \VertexCoordinates{2}{240}
      \RootVertex{0}{0}{#1}
      \DrawVertex{1}{#2}
    \end{tikzpicture}
  }
}

\newcommand{\VertexMM}[3]{
  \vc{
    \begin{tikzpicture}[flag]
      \FlagInit
      \VertexCoordinates{3}{270}
      \RootVertex{0}{0}{#1}
      \DrawVertex{1}{#2}
      \DrawVertex{2}{#3}
      \Connect{0}{1}{black}
      \Connect{0}{2}{black}
      \Connect{1}{2}{black}
    \end{tikzpicture}
  }
}

\newcommand{\VertexNM}[3]{
  \vc{
    \begin{tikzpicture}[flag]
      \FlagInit
      \VertexCoordinates{3}{270}
      \RootVertex{0}{0}{#1}
      \DrawVertex{1}{#2}
      \DrawVertex{2}{#3}
      \Connect{0}{1}{black}
      \Connect{1}{2}{black}
    \end{tikzpicture}
  }
}

\newcommand{\VertexNN}[3]{
  \vc{
    \begin{tikzpicture}[flag]
      \FlagInit
      \VertexCoordinates{3}{270}
      \RootVertex{0}{0}{#1}
      \DrawVertex{1}{#2}
      \DrawVertex{2}{#3}
      \Connect{1}{2}{black}
    \end{tikzpicture}
  }
}

\newcommand{\Edge}[2]{
  \vc{
    \begin{tikzpicture}[flag]
      \FlagInit
      \VertexCoordinates{2}{60}
      \DrawVertex{0}{#1}
      \DrawVertex{1}{#2}
      \Connect{0}{1}{black}
    \end{tikzpicture}
  }
}

\newcommand{\AEdge}[2]{
  \vc{
    \begin{tikzpicture}[flag]
      \FlagInit
      \VertexCoordinates{2}{60}
      \DrawVertex{0}{#1}
      \DrawVertex{1}{#2}
      \ConnectAny{0}{1}{black}
    \end{tikzpicture}
  }
}

\newcommand{\OEdge}[2]{
  \vc{
    \begin{tikzpicture}[flag]
      \FlagInit
      \VertexCoordinates{2}{60}
      \DrawVertex{0}{#1}
      \DrawVertex{1}{#2}
      \ConnectOriented{1}{0}{black}
    \end{tikzpicture}
  }
}

\newcommand{\RootedEdge}[2]{
  \vc{
    \begin{tikzpicture}[flag]
      \FlagInit
      \VertexCoordinates{2}{180}
      \RootVertex{0}{0}{#1}
      \RootVertex{1}{1}{#2}
      \Connect{0}{1}{black}
    \end{tikzpicture}
  }
}

\newcommand{\Triangle}[3]{
  \vc{
    \begin{tikzpicture}[flag]
      \FlagInit
      \VertexCoordinates{3}{210}
      \DrawVertex{0}{#1}
      \DrawVertex{1}{#2}
      \DrawVertex{2}{#3}
      \Connect{0}{1}{black}
      \Connect{1}{2}{black}
      \Connect{0}{2}{black}
    \end{tikzpicture}
  }
}

\newcommand{\OTriangle}[3]{
  \vc{
    \begin{tikzpicture}[flag]
      \FlagInit
      \VertexCoordinates{3}{210}
      \DrawVertex{0}{#1}
      \DrawVertex{1}{#2}
      \DrawVertex{2}{#3}
      \ConnectOriented{0}{1}{black}
      \Connect{1}{2}{black}
      \Connect{0}{2}{black}
    \end{tikzpicture}
  }
}

\newcommand{\TriangleLM}[2]{
  \vc{
    \begin{tikzpicture}[flag]
      \FlagInit
      \VertexCoordinates{6}{210}
      \RootVertex{0}{0}{#1}
      \RootVertex{2}{1}{#1}
      \RootVertex{4}{2}{#1}
      \DrawVertex{5}{#2}
      \Connect{0}{2}{black}
      \Connect{0}{4}{black}
      \Connect{2}{4}{black}
      \Connect{0}{5}{black}
      \Connect{4}{5}{black}
    \end{tikzpicture}
  }
}

\newcommand{\TriangleMR}[2]{
  \vc{
    \begin{tikzpicture}[flag]
      \FlagInit
      \VertexCoordinates{6}{210}
      \RootVertex{0}{0}{#1}
      \RootVertex{2}{1}{#1}
      \RootVertex{4}{2}{#1}
      \DrawVertex{3}{#2}
      \Connect{0}{2}{black}
      \Connect{0}{4}{black}
      \Connect{2}{4}{black}
      \Connect{2}{3}{black}
      \Connect{3}{4}{black}
    \end{tikzpicture}
  }
}

\newcommand{\TriangleLR}[2]{
  \vc{
    \begin{tikzpicture}[flag]
      \FlagInit
      \VertexCoordinates{6}{210}
      \RootVertex{0}{0}{#1}
      \RootVertex{2}{1}{#1}
      \RootVertex{4}{2}{#1}
      \DrawVertex{1}{#2}
      \Connect{0}{2}{black}
      \Connect{0}{4}{black}
      \Connect{2}{4}{black}
      \Connect{0}{1}{black}
      \Connect{1}{2}{black}
    \end{tikzpicture}
  }
}

\newcommand{\EdgeL}[3]{
  \vc{
    \begin{tikzpicture}[flag]
      \FlagInit
      \VertexCoordinates{3}{210}
      \RootVertex{0}{0}{#1}
      \RootVertex{1}{1}{#2}
      \DrawVertex{2}{#3}
      \Connect{0}{1}{black}
      \Connect{0}{2}{black}
    \end{tikzpicture}
  }
}

\newcommand{\OEdgeL}[3]{
  \vc{
    \begin{tikzpicture}[flag]
      \FlagInit
      \VertexCoordinates{3}{210}
      \RootVertex{0}{0}{#1}
      \RootVertex{1}{1}{#2}
      \DrawVertex{2}{#3}
      \ConnectOriented{0}{1}{black}
      \ConnectAny{0}{2}{black}
    \end{tikzpicture}
  }
}

\newcommand{\EdgeR}[3]{
  \vc{
    \begin{tikzpicture}[flag]
      \FlagInit
      \VertexCoordinates{3}{210}
      \RootVertex{0}{0}{#1}
      \RootVertex{1}{1}{#2}
      \DrawVertex{2}{#3}
      \Connect{0}{1}{black}
      \Connect{1}{2}{black}
    \end{tikzpicture}
  }
}

\newcommand{\OEdgeR}[3]{
  \vc{
    \begin{tikzpicture}[flag]
      \FlagInit
      \VertexCoordinates{3}{210}
      \RootVertex{0}{0}{#1}
      \RootVertex{1}{1}{#2}
      \DrawVertex{2}{#3}
      \ConnectOriented{0}{1}{black}
      \ConnectAny{1}{2}{black}
    \end{tikzpicture}
  }
}

\newcommand{\EdgeN}[3]{
  \vc{
    \begin{tikzpicture}[flag]
      \FlagInit
      \VertexCoordinates{3}{210}
      \RootVertex{0}{0}{#1}
      \RootVertex{1}{1}{#2}
      \DrawVertex{2}{#3}
      \Connect{0}{1}{black}
    \end{tikzpicture}
  }
}

\newcommand{\OEdgeN}[3]{
  \vc{
    \begin{tikzpicture}[flag]
      \FlagInit
      \VertexCoordinates{3}{210}
      \RootVertex{0}{0}{#1}
      \RootVertex{1}{1}{#2}
      \DrawVertex{2}{#3}
      \ConnectOriented{0}{1}{black}
    \end{tikzpicture}
  }
}

\newcommand{\EdgeB}[3]{
  \vc{
    \begin{tikzpicture}[flag]
      \FlagInit
      \VertexCoordinates{3}{210}
      \RootVertex{0}{0}{#1}
      \RootVertex{1}{1}{#2}
      \DrawVertex{2}{#3}
      \Connect{0}{1}{black}
      \Connect{0}{2}{black}
      \Connect{1}{2}{black}
    \end{tikzpicture}
  }
}

\newcommand{\EdgeBOA}[3]{
  \vc{
    \begin{tikzpicture}[flag]
      \FlagInit
      \VertexCoordinates{3}{210}
      \RootVertex{0}{0}{#1}
      \RootVertex{1}{1}{#2}
      \DrawVertex{2}{#3}
      \Connect{0}{1}{black}
      \ConnectOriented{0}{2}{black}
      \ConnectAny{1}{2}{black}
    \end{tikzpicture}
  }
}

\newcommand{\EdgeBOO}[3]{
  \vc{
    \begin{tikzpicture}[flag]
      \FlagInit
      \VertexCoordinates{3}{210}
      \RootVertex{0}{0}{#1}
      \RootVertex{1}{1}{#2}
      \DrawVertex{2}{#3}
      \Connect{0}{1}{black}
      \ConnectOriented{0}{2}{black}
      \ConnectOriented{1}{2}{black}
    \end{tikzpicture}
  }
}

\newcommand{\EdgeBOo}[3]{
  \vc{
    \begin{tikzpicture}[flag]
      \FlagInit
      \VertexCoordinates{3}{210}
      \RootVertex{0}{0}{#1}
      \RootVertex{1}{1}{#2}
      \DrawVertex{2}{#3}
      \Connect{0}{1}{black}
      \ConnectOriented{0}{2}{black}
      \ConnectOriented{2}{1}{black}
    \end{tikzpicture}
  }
}

\newcommand{\EdgeBoO}[3]{
  \vc{
    \begin{tikzpicture}[flag]
      \FlagInit
      \VertexCoordinates{3}{210}
      \RootVertex{0}{0}{#1}
      \RootVertex{1}{1}{#2}
      \DrawVertex{2}{#3}
      \Connect{0}{1}{black}
      \ConnectOriented{2}{0}{black}
      \ConnectOriented{1}{2}{black}
    \end{tikzpicture}
  }
}

\newcommand{\EdgeBoo}[3]{
  \vc{
    \begin{tikzpicture}[flag]
      \FlagInit
      \VertexCoordinates{3}{210}
      \RootVertex{0}{0}{#1}
      \RootVertex{1}{1}{#2}
      \DrawVertex{2}{#3}
      \Connect{0}{1}{black}
      \ConnectOriented{2}{0}{black}
      \ConnectOriented{2}{1}{black}
    \end{tikzpicture}
  }
}

\newcommand{\EdgeBON}[3]{
  \vc{
    \begin{tikzpicture}[flag]
      \FlagInit
      \VertexCoordinates{3}{210}
      \RootVertex{0}{0}{#1}
      \RootVertex{1}{1}{#2}
      \DrawVertex{2}{#3}
      \Connect{0}{1}{black}
      \ConnectOriented{0}{2}{black}
      \Connect{1}{2}{black}
    \end{tikzpicture}
  }
}

\newcommand{\EdgeBoN}[3]{
  \vc{
    \begin{tikzpicture}[flag]
      \FlagInit
      \VertexCoordinates{3}{210}
      \RootVertex{0}{0}{#1}
      \RootVertex{1}{1}{#2}
      \DrawVertex{2}{#3}
      \Connect{0}{1}{black}
      \ConnectOriented{2}{0}{black}
      \Connect{1}{2}{black}
    \end{tikzpicture}
  }
}

\newcommand{\OEdgeB}[3]{
  \vc{
    \begin{tikzpicture}[flag]
      \FlagInit
      \VertexCoordinates{3}{210}
      \RootVertex{0}{0}{#1}
      \RootVertex{1}{1}{#2}
      \DrawVertex{2}{#3}
      \ConnectOriented{0}{1}{black}
      \ConnectAny{0}{2}{black}
      \ConnectAny{1}{2}{black}
    \end{tikzpicture}
  }
}

\newcommand{\EdgeNaN}[4]{
  \vc{
    \begin{tikzpicture}[flag]
      \FlagInit
      \VertexCoordinates{4}{225}
      \RootVertex{0}{0}{#1}
      \RootVertex{1}{1}{#2}
      \DrawVertex{2}{#3}
      \DrawVertex{3}{#4}
      \Connect{0}{1}{black}
      \Connect{2}{3}{black}
    \end{tikzpicture}
  }
}

\newcommand{\OEdgeNaN}[4]{
  \vc{
    \begin{tikzpicture}[flag]
      \FlagInit
      \VertexCoordinates{4}{225}
      \RootVertex{0}{0}{#1}
      \RootVertex{1}{1}{#2}
      \DrawVertex{2}{#3}
      \DrawVertex{3}{#4}
      \ConnectOriented{0}{1}{black}
      \ConnectAny{2}{3}{black}
    \end{tikzpicture}
  }
}

\newcommand{\EdgeLN}[4]{
  \vc{
    \begin{tikzpicture}[flag]
      \FlagInit
      \VertexCoordinates{4}{225}
      \RootVertex{0}{0}{#1}
      \RootVertex{1}{1}{#2}
      \DrawVertex{2}{#3}
      \DrawVertex{3}{#4}
      \Connect{0}{1}{black}
      \Connect{0}{3}{black}
      \Connect{2}{3}{black}
    \end{tikzpicture}
  }
}

\newcommand{\OEdgeLN}[4]{
  \vc{
    \begin{tikzpicture}[flag]
      \FlagInit
      \VertexCoordinates{4}{225}
      \RootVertex{0}{0}{#1}
      \RootVertex{1}{1}{#2}
      \DrawVertex{2}{#3}
      \DrawVertex{3}{#4}
      \ConnectOriented{0}{1}{black}
      \ConnectAny{0}{3}{black}
      \ConnectAny{2}{3}{black}
    \end{tikzpicture}
  }
}

\newcommand{\EdgeRN}[4]{
  \vc{
    \begin{tikzpicture}[flag]
      \FlagInit
      \VertexCoordinates{4}{225}
      \RootVertex{0}{0}{#1}
      \RootVertex{1}{1}{#2}
      \DrawVertex{2}{#3}
      \DrawVertex{3}{#4}
      \Connect{0}{1}{black}
      \Connect{1}{2}{black}
      \Connect{2}{3}{black}
    \end{tikzpicture}
  }
}

\newcommand{\OEdgeRN}[4]{
  \vc{
    \begin{tikzpicture}[flag]
      \FlagInit
      \VertexCoordinates{4}{225}
      \RootVertex{0}{0}{#1}
      \RootVertex{1}{1}{#2}
      \DrawVertex{2}{#3}
      \DrawVertex{3}{#4}
      \ConnectOriented{0}{1}{black}
      \ConnectAny{1}{2}{black}
      \ConnectAny{2}{3}{black}
    \end{tikzpicture}
  }
}

\newcommand{\EdgeBN}[4]{
  \vc{
    \begin{tikzpicture}[flag]
      \FlagInit
      \VertexCoordinates{4}{225}
      \RootVertex{0}{0}{#1}
      \RootVertex{1}{1}{#2}
      \DrawVertex{2}{#3}
      \DrawVertex{3}{#4}
      \Connect{0}{1}{black}
      \Connect{0}{3}{black}
      \Connect{1}{3}{black}
      \Connect{2}{3}{black}
    \end{tikzpicture}
  }
}

\newcommand{\OEdgeBN}[4]{
  \vc{
    \begin{tikzpicture}[flag]
      \FlagInit
      \VertexCoordinates{4}{225}
      \RootVertex{0}{0}{#1}
      \RootVertex{1}{1}{#2}
      \DrawVertex{2}{#3}
      \DrawVertex{3}{#4}
      \ConnectOriented{0}{1}{black}
      \ConnectAny{0}{3}{black}
      \ConnectAny{1}{3}{black}
      \ConnectAny{2}{3}{black}
    \end{tikzpicture}
  }
}

\newcommand{\EdgeLaR}[4]{
  \vc{
    \begin{tikzpicture}[flag]
      \FlagInit
      \VertexCoordinates{4}{225}
      \RootVertex{0}{0}{#1}
      \RootVertex{1}{1}{#2}
      \DrawVertex{2}{#3}
      \DrawVertex{3}{#4}
      \Connect{0}{1}{black}
      \Connect{0}{3}{black}
      \Connect{1}{2}{black}
      \Connect{2}{3}{black}
    \end{tikzpicture}
  }
}

\newcommand{\EdgeBL}[4]{
  \vc{
    \begin{tikzpicture}[flag]
      \FlagInit
      \VertexCoordinates{4}{225}
      \RootVertex{0}{0}{#1}
      \RootVertex{1}{1}{#2}
      \DrawVertex{2}{#3}
      \DrawVertex{3}{#4}
      \Connect{0}{1}{black}
      \Connect{0}{2}{black}
      \Connect{1}{2}{black}
      \Connect{0}{3}{black}
      \Connect{2}{3}{black}
    \end{tikzpicture}
  }
}

\newcommand{\OEdgeBL}[4]{
  \vc{
    \begin{tikzpicture}[flag]
      \FlagInit
      \VertexCoordinates{4}{225}
      \RootVertex{0}{0}{#1}
      \RootVertex{1}{1}{#2}
      \DrawVertex{2}{#3}
      \DrawVertex{3}{#4}
      \ConnectOriented{0}{1}{black}
      \ConnectAny{0}{2}{black}
      \ConnectAny{1}{2}{black}
      \ConnectAny{0}{3}{black}
      \ConnectAny{2}{3}{black}
    \end{tikzpicture}
  }
}

\newcommand{\EdgeBR}[4]{
  \vc{
    \begin{tikzpicture}[flag]
      \FlagInit
      \VertexCoordinates{4}{225}
      \RootVertex{0}{0}{#1}
      \RootVertex{1}{1}{#2}
      \DrawVertex{2}{#3}
      \DrawVertex{3}{#4}
      \Connect{0}{1}{black}
      \Connect{1}{2}{black}
      \Connect{0}{3}{black}
      \Connect{1}{3}{black}
      \Connect{2}{3}{black}
    \end{tikzpicture}
  }
}

\newcommand{\OEdgeBR}[4]{
  \vc{
    \begin{tikzpicture}[flag]
      \FlagInit
      \VertexCoordinates{4}{225}
      \RootVertex{0}{0}{#1}
      \RootVertex{1}{1}{#2}
      \DrawVertex{2}{#3}
      \DrawVertex{3}{#4}
      \ConnectOriented{0}{1}{black}
      \ConnectAny{1}{2}{black}
      \ConnectAny{0}{3}{black}
      \ConnectAny{1}{3}{black}
      \ConnectAny{2}{3}{black}
    \end{tikzpicture}
  }
}

\newcommand{\EdgeLaL}[4]{
  \vc{
    \begin{tikzpicture}[flag]
      \FlagInit
      \VertexCoordinates{4}{225}
      \RootVertex{0}{0}{#1}
      \RootVertex{1}{1}{#2}
      \DrawVertex{2}{#3}
      \DrawVertex{3}{#4}
      \Connect{0}{1}{black}
      \Connect{0}{2}{black}
      \Connect{0}{3}{black}
      \Connect{2}{3}{black}
    \end{tikzpicture}
  }
}

\newcommand{\OEdgeLaL}[4]{
  \vc{
    \begin{tikzpicture}[flag]
      \FlagInit
      \VertexCoordinates{4}{225}
      \RootVertex{0}{0}{#1}
      \RootVertex{1}{1}{#2}
      \DrawVertex{2}{#3}
      \DrawVertex{3}{#4}
      \ConnectOriented{0}{1}{black}
      \ConnectAny{0}{2}{black}
      \ConnectAny{0}{3}{black}
      \ConnectAny{2}{3}{black}
    \end{tikzpicture}
  }
}

\newcommand{\EdgeRaR}[4]{
  \vc{
    \begin{tikzpicture}[flag]
      \FlagInit
      \VertexCoordinates{4}{225}
      \RootVertex{0}{0}{#1}
      \RootVertex{1}{1}{#2}
      \DrawVertex{2}{#3}
      \DrawVertex{3}{#4}
      \Connect{0}{1}{black}
      \Connect{1}{2}{black}
      \Connect{1}{3}{black}
      \Connect{2}{3}{black}
    \end{tikzpicture}
  }
}

\newcommand{\OEdgeRaR}[4]{
  \vc{
    \begin{tikzpicture}[flag]
      \FlagInit
      \VertexCoordinates{4}{225}
      \RootVertex{0}{0}{#1}
      \RootVertex{1}{1}{#2}
      \DrawVertex{2}{#3}
      \DrawVertex{3}{#4}
      \ConnectOriented{0}{1}{black}
      \ConnectAny{1}{2}{black}
      \ConnectAny{1}{3}{black}
      \ConnectAny{2}{3}{black}
    \end{tikzpicture}
  }
}

\newcommand{\CherryN}[2]{
  \vc{
    \begin{tikzpicture}[flag]
      \FlagInit
      \VertexCoordinates{4}{180}
      \RootVertex{0}{1}{#1}
      \RootVertex{1}{0}{#1}
      \RootVertex{2}{2}{#1}
      \DrawVertex{3}{#2}
      \Connect{0}{1}{black}
      \Connect{1}{2}{black}
    \end{tikzpicture}
  }
}

\newcommand{\CherryM}[2]{
  \vc{
    \begin{tikzpicture}[flag]
      \FlagInit
      \VertexCoordinates{4}{180}
      \RootVertex{0}{1}{#1}
      \RootVertex{1}{0}{#1}
      \RootVertex{2}{2}{#1}
      \DrawVertex{3}{#2}
      \Connect{0}{1}{black}
      \Connect{1}{2}{black}
      \Connect{1}{3}{black}
    \end{tikzpicture}
  }
}

\newcommand{\CherryR}[2]{
  \vc{
    \begin{tikzpicture}[flag]
      \FlagInit
      \VertexCoordinates{4}{180}
      \RootVertex{0}{1}{#1}
      \RootVertex{1}{0}{#1}
      \RootVertex{2}{2}{#1}
      \DrawVertex{3}{#2}
      \Connect{0}{1}{black}
      \Connect{1}{2}{black}
      \Connect{2}{3}{black}
    \end{tikzpicture}
  }
}

\newcommand{\CherryL}[2]{
  \vc{
    \begin{tikzpicture}[flag]
      \FlagInit
      \VertexCoordinates{4}{180}
      \RootVertex{0}{1}{#1}
      \RootVertex{1}{0}{#1}
      \RootVertex{2}{2}{#1}
      \DrawVertex{3}{#2}
      \Connect{0}{1}{black}
      \Connect{1}{2}{black}
      \Connect{0}{3}{black}
    \end{tikzpicture}
  }
}

\newcommand{\CherryLR}[2]{
  \vc{
    \begin{tikzpicture}[flag]
      \FlagInit
      \VertexCoordinates{4}{180}
      \RootVertex{0}{1}{#1}
      \RootVertex{1}{0}{#1}
      \RootVertex{2}{2}{#1}
      \DrawVertex{3}{#2}
      \Connect{0}{1}{black}
      \Connect{1}{2}{black}
      \Connect{0}{3}{black}
      \Connect{2}{3}{black}
    \end{tikzpicture}
  }
}

\newcommand{\CherryMR}[2]{
  \vc{
    \begin{tikzpicture}[flag]
      \FlagInit
      \VertexCoordinates{4}{180}
      \RootVertex{0}{1}{#1}
      \RootVertex{1}{0}{#1}
      \RootVertex{2}{2}{#1}
      \DrawVertex{3}{#2}
      \Connect{0}{1}{black}
      \Connect{1}{2}{black}
      \Connect{1}{3}{black}
      \Connect{2}{3}{black}
    \end{tikzpicture}
  }
}

\newcommand{\CherryLM}[2]{
  \vc{
    \begin{tikzpicture}[flag]
      \FlagInit
      \VertexCoordinates{4}{180}
      \RootVertex{0}{1}{#1}
      \RootVertex{1}{0}{#1}
      \RootVertex{2}{2}{#1}
      \DrawVertex{3}{#2}
      \Connect{0}{1}{black}
      \Connect{1}{2}{black}
      \Connect{0}{3}{black}
      \Connect{1}{3}{black}
    \end{tikzpicture}
  }
}

\newcommand{\CherryA}[2]{
  \vc{
    \begin{tikzpicture}[flag]
      \FlagInit
      \VertexCoordinates{4}{180}
      \RootVertex{0}{1}{#1}
      \RootVertex{1}{0}{#1}
      \RootVertex{2}{2}{#1}
      \DrawVertex{3}{#2}
      \Connect{0}{1}{black}
      \Connect{1}{2}{black}
      \Connect{0}{3}{black}
      \Connect{1}{3}{black}
      \Connect{2}{3}{black}
    \end{tikzpicture}
  }
}

\newcommand{\CherryLaL}[2]{
  \vc{
    \begin{tikzpicture}[flag]
      \FlagInit
      \VertexCoordinates{5}{198}
      \RootVertex{0}{1}{#1}
      \RootVertex{1}{0}{#1}
      \RootVertex{2}{2}{#1}
      \DrawVertex{3}{#2}
      \DrawVertex{4}{#2}
      \Connect{0}{1}{black}
      \Connect{1}{2}{black}
      \Connect{0}{3}{black}
      \Connect{0}{4}{black}
      \Connect{3}{4}{black}
    \end{tikzpicture}
  }
}

\newcommand{\CherryRaR}[2]{
  \vc{
    \begin{tikzpicture}[flag]
      \FlagInit
      \VertexCoordinates{5}{198}
      \RootVertex{0}{1}{#1}
      \RootVertex{1}{0}{#1}
      \RootVertex{2}{2}{#1}
      \DrawVertex{3}{#2}
      \DrawVertex{4}{#2}
      \Connect{0}{1}{black}
      \Connect{1}{2}{black}
      \Connect{2}{3}{black}
      \Connect{2}{4}{black}
      \Connect{3}{4}{black}
    \end{tikzpicture}
  }
}

\newcommand{\CherryLaR}[2]{
  \vc{
    \begin{tikzpicture}[flag]
      \FlagInit
      \VertexCoordinates{5}{198}
      \RootVertex{0}{1}{#1}
      \RootVertex{1}{0}{#1}
      \RootVertex{2}{2}{#1}
      \DrawVertex{3}{#2}
      \DrawVertex{4}{#2}
      \Connect{0}{1}{black}
      \Connect{1}{2}{black}
      \Connect{2}{3}{black}
      \Connect{0}{4}{black}
      \Connect{3}{4}{black}
    \end{tikzpicture}
  }
}

\newcommand{\CherryLRaR}[2]{
  \vc{
    \begin{tikzpicture}[flag]
      \FlagInit
      \VertexCoordinates{5}{198}
      \RootVertex{0}{1}{#1}
      \RootVertex{1}{0}{#1}
      \RootVertex{2}{2}{#1}
      \DrawVertex{3}{#2}
      \DrawVertex{4}{#2}
      \Connect{0}{1}{black}
      \Connect{1}{2}{black}
      \Connect{2}{3}{black}
      \Connect{0}{4}{black}
      \Connect{2}{4}{black}
      \Connect{3}{4}{black}
    \end{tikzpicture}
  }
}

\newcommand{\CherryLaLR}[2]{
  \vc{
    \begin{tikzpicture}[flag]
      \FlagInit
      \VertexCoordinates{5}{198}
      \RootVertex{0}{1}{#1}
      \RootVertex{1}{0}{#1}
      \RootVertex{2}{2}{#1}
      \DrawVertex{3}{#2}
      \DrawVertex{4}{#2}
      \Connect{0}{1}{black}
      \Connect{1}{2}{black}
      \Connect{0}{3}{black}
      \Connect{2}{3}{black}
      \Connect{0}{4}{black}
      \Connect{3}{4}{black}
    \end{tikzpicture}
  }
}

\newcommand{\CherryLRaLR}[2]{
  \vc{
    \begin{tikzpicture}[flag]
      \FlagInit
      \VertexCoordinates{5}{198}
      \RootVertex{0}{1}{#1}
      \RootVertex{1}{0}{#1}
      \RootVertex{2}{2}{#1}
      \DrawVertex{3}{#2}
      \DrawVertex{4}{#2}
      \Connect{0}{1}{black}
      \Connect{1}{2}{black}
      \Connect{0}{3}{black}
      \Connect{2}{3}{black}
      \Connect{0}{4}{black}
      \Connect{2}{4}{black}
      \Connect{3}{4}{black}
    \end{tikzpicture}
  }
}

\newcommand{\CherryNN}[2]{
  \vc{
    \begin{tikzpicture}[flag]
      \FlagInit
      \VertexCoordinates{5}{198}
      \RootVertex{0}{1}{#1}
      \RootVertex{1}{0}{#1}
      \RootVertex{2}{2}{#1}
      \DrawVertex{3}{#2}
      \DrawVertex{4}{#2}
      \Connect{0}{1}{black}
      \Connect{1}{2}{black}
      \Connect{3}{4}{black}
    \end{tikzpicture}
  }
}

\newcommand{\CherryNA}[2]{
  \vc{
    \begin{tikzpicture}[flag]
      \FlagInit
      \VertexCoordinates{5}{198}
      \RootVertex{0}{1}{#1}
      \RootVertex{1}{0}{#1}
      \RootVertex{2}{2}{#1}
      \DrawVertex{3}{#2}
      \DrawVertex{4}{#2}
      \Connect{0}{1}{black}
      \Connect{1}{2}{black}
      \Connect{0}{3}{black}
      \Connect{1}{3}{black}
      \Connect{2}{3}{black}
      \Connect{3}{4}{black}
    \end{tikzpicture}
  }
}

\newcommand{\CherryNM}[2]{
  \vc{
    \begin{tikzpicture}[flag]
      \FlagInit
      \VertexCoordinates{5}{198}
      \RootVertex{0}{1}{#1}
      \RootVertex{1}{0}{#1}
      \RootVertex{2}{2}{#1}
      \DrawVertex{3}{#2}
      \DrawVertex{4}{#2}
      \Connect{0}{1}{black}
      \Connect{1}{2}{black}
      \Connect{1}{3}{black}
      \Connect{3}{4}{black}
    \end{tikzpicture}
  }
}

\newcommand{\CherryMM}[2]{
  \vc{
    \begin{tikzpicture}[flag]
      \FlagInit
      \VertexCoordinates{5}{198}
      \RootVertex{0}{1}{#1}
      \RootVertex{1}{0}{#1}
      \RootVertex{2}{2}{#1}
      \DrawVertex{3}{#2}
      \DrawVertex{4}{#2}
      \Connect{0}{1}{black}
      \Connect{1}{2}{black}
      \Connect{1}{3}{black}
      \Connect{1}{4}{black}
      \Connect{3}{4}{black}
    \end{tikzpicture}
  }
}

\newcommand{\CherryMRaR}[2]{
  \vc{
    \begin{tikzpicture}[flag]
      \FlagInit
      \VertexCoordinates{5}{198}
      \RootVertex{0}{1}{#1}
      \RootVertex{1}{0}{#1}
      \RootVertex{2}{2}{#1}
      \DrawVertex{3}{#2}
      \DrawVertex{4}{#2}
      \Connect{0}{1}{black}
      \Connect{1}{2}{black}
      \Connect{2}{3}{black}
      \Connect{1}{4}{black}
      \Connect{2}{4}{black}
      \Connect{3}{4}{black}
    \end{tikzpicture}
  }
}

\newcommand{\CherryLMaR}[2]{
  \vc{
    \begin{tikzpicture}[flag]
      \FlagInit
      \VertexCoordinates{5}{198}
      \RootVertex{0}{1}{#1}
      \RootVertex{1}{0}{#1}
      \RootVertex{2}{2}{#1}
      \DrawVertex{3}{#2}
      \DrawVertex{4}{#2}
      \Connect{0}{1}{black}
      \Connect{1}{2}{black}
      \Connect{2}{3}{black}
      \Connect{0}{4}{black}
      \Connect{1}{4}{black}
      \Connect{3}{4}{black}
    \end{tikzpicture}
  }
}

\newcommand{\CherryAaR}[2]{
  \vc{
    \begin{tikzpicture}[flag]
      \FlagInit
      \VertexCoordinates{5}{198}
      \RootVertex{0}{1}{#1}
      \RootVertex{1}{0}{#1}
      \RootVertex{2}{2}{#1}
      \DrawVertex{3}{#2}
      \DrawVertex{4}{#2}
      \Connect{0}{1}{black}
      \Connect{1}{2}{black}
      \Connect{2}{3}{black}
      \Connect{0}{4}{black}
      \Connect{1}{4}{black}
      \Connect{2}{4}{black}
      \Connect{3}{4}{black}
    \end{tikzpicture}
  }
}

\newcommand{\CherryLaMR}[2]{
  \vc{
    \begin{tikzpicture}[flag]
      \FlagInit
      \VertexCoordinates{5}{198}
      \RootVertex{0}{1}{#1}
      \RootVertex{1}{0}{#1}
      \RootVertex{2}{2}{#1}
      \DrawVertex{3}{#2}
      \DrawVertex{4}{#2}
      \Connect{0}{1}{black}
      \Connect{1}{2}{black}
      \Connect{1}{3}{black}
      \Connect{2}{3}{black}
      \Connect{0}{4}{black}
      \Connect{3}{4}{black}
    \end{tikzpicture}
  }
}

\newcommand{\CherryLRaMR}[2]{
  \vc{
    \begin{tikzpicture}[flag]
      \FlagInit
      \VertexCoordinates{5}{198}
      \RootVertex{0}{1}{#1}
      \RootVertex{1}{0}{#1}
      \RootVertex{2}{2}{#1}
      \DrawVertex{3}{#2}
      \DrawVertex{4}{#2}
      \Connect{0}{1}{black}
      \Connect{1}{2}{black}
      \Connect{1}{3}{black}
      \Connect{2}{3}{black}
      \Connect{0}{4}{black}
      \Connect{2}{4}{black}
      \Connect{3}{4}{black}
    \end{tikzpicture}
  }
}

\newcommand{\CherryLaLM}[2]{
  \vc{
    \begin{tikzpicture}[flag]
      \FlagInit
      \VertexCoordinates{5}{198}
      \RootVertex{0}{1}{#1}
      \RootVertex{1}{0}{#1}
      \RootVertex{2}{2}{#1}
      \DrawVertex{3}{#2}
      \DrawVertex{4}{#2}
      \Connect{0}{1}{black}
      \Connect{1}{2}{black}
      \Connect{0}{3}{black}
      \Connect{1}{3}{black}
      \Connect{0}{4}{black}
      \Connect{3}{4}{black}
    \end{tikzpicture}
  }
}

\newcommand{\CherryLaA}[2]{
  \vc{
    \begin{tikzpicture}[flag]
      \FlagInit
      \VertexCoordinates{5}{198}
      \RootVertex{0}{1}{#1}
      \RootVertex{1}{0}{#1}
      \RootVertex{2}{2}{#1}
      \DrawVertex{3}{#2}
      \DrawVertex{4}{#2}
      \Connect{0}{1}{black}
      \Connect{1}{2}{black}
      \Connect{0}{3}{black}
      \Connect{1}{3}{black}
      \Connect{2}{3}{black}
      \Connect{0}{4}{black}
      \Connect{3}{4}{black}
    \end{tikzpicture}
  }
}

\newcommand{\CherryLMaLR}[2]{
  \vc{
    \begin{tikzpicture}[flag]
      \FlagInit
      \VertexCoordinates{5}{198}
      \RootVertex{0}{1}{#1}
      \RootVertex{1}{0}{#1}
      \RootVertex{2}{2}{#1}
      \DrawVertex{3}{#2}
      \DrawVertex{4}{#2}
      \Connect{0}{1}{black}
      \Connect{1}{2}{black}
      \Connect{0}{3}{black}
      \Connect{2}{3}{black}
      \Connect{0}{4}{black}
      \Connect{1}{4}{black}
      \Connect{3}{4}{black}
    \end{tikzpicture}
  }
}

\newcommand{\CherryLRaA}[2]{
  \vc{
    \begin{tikzpicture}[flag]
      \FlagInit
      \VertexCoordinates{5}{198}
      \RootVertex{0}{1}{#1}
      \RootVertex{1}{0}{#1}
      \RootVertex{2}{2}{#1}
      \DrawVertex{3}{#2}
      \DrawVertex{4}{#2}
      \Connect{0}{1}{black}
      \Connect{1}{2}{black}
      \Connect{0}{3}{black}
      \Connect{1}{3}{black}
      \Connect{2}{3}{black}
      \Connect{0}{4}{black}
      \Connect{2}{4}{black}
      \Connect{3}{4}{black}
    \end{tikzpicture}
  }
}

\newcommand{\CherryNaMR}[2]{
  \vc{
    \begin{tikzpicture}[flag]
      \FlagInit
      \VertexCoordinates{5}{198}
      \RootVertex{0}{1}{#1}
      \RootVertex{1}{0}{#1}
      \RootVertex{2}{2}{#1}
      \DrawVertex{3}{#2}
      \DrawVertex{4}{#2}
      \Connect{0}{1}{black}
      \Connect{1}{2}{black}
      \Connect{1}{3}{black}
      \Connect{2}{3}{black}
      \Connect{3}{4}{black}
    \end{tikzpicture}
  }
}

\newcommand{\CherryLMaN}[2]{
  \vc{
    \begin{tikzpicture}[flag]
      \FlagInit
      \VertexCoordinates{5}{198}
      \RootVertex{0}{1}{#1}
      \RootVertex{1}{0}{#1}
      \RootVertex{2}{2}{#1}
      \DrawVertex{3}{#2}
      \DrawVertex{4}{#2}
      \Connect{0}{1}{black}
      \Connect{1}{2}{black}
      \Connect{0}{4}{black}
      \Connect{1}{4}{black}
      \Connect{3}{4}{black}
    \end{tikzpicture}
  }
}

\newcommand{\CherryNaA}[2]{
  \vc{
    \begin{tikzpicture}[flag]
      \FlagInit
      \VertexCoordinates{5}{198}
      \RootVertex{0}{1}{#1}
      \RootVertex{1}{0}{#1}
      \RootVertex{2}{2}{#1}
      \DrawVertex{3}{#2}
      \DrawVertex{4}{#2}
      \Connect{0}{1}{black}
      \Connect{1}{2}{black}
      \Connect{0}{3}{black}
      \Connect{1}{3}{black}
      \Connect{2}{3}{black}
      \Connect{3}{4}{black}
    \end{tikzpicture}
  }
}

\newcommand{\CherryMaMR}[2]{
  \vc{
    \begin{tikzpicture}[flag]
      \FlagInit
      \VertexCoordinates{5}{198}
      \RootVertex{0}{1}{#1}
      \RootVertex{1}{0}{#1}
      \RootVertex{2}{2}{#1}
      \DrawVertex{3}{#2}
      \DrawVertex{4}{#2}
      \Connect{0}{1}{black}
      \Connect{1}{2}{black}
      \Connect{1}{3}{black}
      \Connect{2}{3}{black}
      \Connect{1}{4}{black}
      \Connect{3}{4}{black}
    \end{tikzpicture}
  }
}

\newcommand{\CherryLMaM}[2]{
  \vc{
    \begin{tikzpicture}[flag]
      \FlagInit
      \VertexCoordinates{5}{198}
      \RootVertex{0}{1}{#1}
      \RootVertex{1}{0}{#1}
      \RootVertex{2}{2}{#1}
      \DrawVertex{3}{#2}
      \DrawVertex{4}{#2}
      \Connect{0}{1}{black}
      \Connect{1}{2}{black}
      \Connect{1}{3}{black}
      \Connect{0}{4}{black}
      \Connect{1}{4}{black}
      \Connect{3}{4}{black}
    \end{tikzpicture}
  }
}

\newcommand{\CherryMaA}[2]{
  \vc{
    \begin{tikzpicture}[flag]
      \FlagInit
      \VertexCoordinates{5}{198}
      \RootVertex{0}{1}{#1}
      \RootVertex{1}{0}{#1}
      \RootVertex{2}{2}{#1}
      \DrawVertex{3}{#2}
      \DrawVertex{4}{#2}
      \Connect{0}{1}{black}
      \Connect{1}{2}{black}
      \Connect{0}{3}{black}
      \Connect{1}{3}{black}
      \Connect{2}{3}{black}
      \Connect{1}{4}{black}
      \Connect{3}{4}{black}
    \end{tikzpicture}
  }
}

\newcommand{\CherryLMaMR}[2]{
  \vc{
    \begin{tikzpicture}[flag]
      \FlagInit
      \VertexCoordinates{5}{198}
      \RootVertex{0}{1}{#1}
      \RootVertex{1}{0}{#1}
      \RootVertex{2}{2}{#1}
      \DrawVertex{3}{#2}
      \DrawVertex{4}{#2}
      \Connect{0}{1}{black}
      \Connect{1}{2}{black}
      \Connect{1}{3}{black}
      \Connect{2}{3}{black}
      \Connect{0}{4}{black}
      \Connect{1}{4}{black}
      \Connect{3}{4}{black}
    \end{tikzpicture}
  }
}

%% file: sections/introduction.tex
The \emph{Sparse half conjecture} of Erd\H{o}s (\cite{Original}) states that in a triangle-free graph on $n$ vertices one can always find a vertex subset of size at least $\frac{n}{2}$ such that at most $\frac{n^2}{50}$ edges have both of their endpoints in it.
This conjecture has received a lot of attention (\cites{tri-dense, tri-sparse, local-density-condition, tri-sparse-first}).
The best bound so far is $\frac{27n^2}{1024}$, and it was achieved by Razborov in \cite{tri-sparse}.

A related conjecture of Erd\H{o}s is that every triangle-free graph on $n$ vertices can be made bipartite by removing at most $\frac{n^2}{25}$ edges.
This question has also been extensively studied over the years (\cites{make-bipartite, make-bipartite-2, sparse-halves-best}), and the
best known bound of $\frac{2n^2}{47}$ was proven in \cite{sparse-halves-best}.

%Let $G$ be a simple graph with vertex set $V(G)$ and edge set $E(G)$.
%By $E(G[A])$ we denote the set of edges with both endpoints in the set $A \subseteq V(G)$,
%while by $e(G[A])$ we denote the number of edges with both endpoints in $A$.
%Analogously, for two disjoint sets $A, B \subseteq V(G)$, $E(A, B)$ is the set of edges with one endpoint in $A$ and the other in $B$ and $e(G[A], G[B])$ is the number of such edges.
%When $G$ is clear from context, we abuse notation and use $E(A)$, $e(A)$, $E(A, B)$, and $e(A, B)$, respectively.
%A vertex $k$-partition of a graph $G$ is a partition of its vertex set into $k$ sets $A_1, \ldots, A_k$.
%We say that a $k$-partition is balanced if $||A_i| - |A_j|| < 1$ for every pair of integers $i, j \leq k$.
%An edge $(u, v) \in E(G)$ is called a class-edge, if there exists an $i$, such that $\{u, v\} \subseteq A_i$.

Recently, Balogh, Clemen and Lidick\'y (\cite{10problems}) have proposed a number of problems on finding partitions with few edges inside the parts.
Many of these focus on balanced $k$-partitions, partly because a good balanced bipartition yields a sparse half.
We recall that a $k$-partition $V(G) = A_1 \cup \ldots \cup A_k$ is balanced if $||A_i| - |A_j|| \leq 1$ for all $i, j \in \{1, \ldots, k\}$. 
Let 
\begin{equation}\label{k-distance-definition}
  D^b_k(G) = \min_{A_1, \ldots, A_k} \sum_{i = 1}^k e(A_i),
\end{equation}
where the minimum is taken over all balanced $k$-partitions and $e(A_i)$ denotes the number of edges with both endpoints in $A_i$.
It turns out that proving tight bounds for small values of $k$ in the class of triangle-free graphs is easier than resolving the original conjectures of Erd\H{o}s.

\begin{theorem}[\cite{10problems}]
  If $G$ is a sufficiently large triangle-free graph on $n$ vertices, then
  $$D^b_2(G) \leq \frac{n^2}{16} \qquad\text{and}\qquad D^b_3(G) \leq \frac{n^2}{36}.$$
  Moreover, $D^b_2(K_{\frac{3n}{4},\frac{n}{4}}) = \frac{n^2}{16}$ and $D^b_3(K_{\frac{n}{2},\frac{n}{2}}) = \frac{n^2}{36}$.
\end{theorem}

The authors of \cite{10problems} also considered such problems in the class of $K_4$-free graphs and formulated a conjecture on the optimal bound for $D_2^b$.
This conjecture was recently proved by Balogh and the authors of this paper. 

\begin{theorem}[\cite{balanced-bipartite-distance-k4-free}]
  If $G$ is a $K_4$-free graph on $n$ vertices, then
  $$D^b_2(G) \leq \frac{n^2}{9}.$$
  Moreover, $D^b_2(K_{\frac{n}{3},\frac{n}{3},\frac{n}{3}}) = \frac{n^2}{9}$.
\end{theorem}

A natural way of generalizing \eqref{k-distance-definition} is by considering the $\ell_p$-norm of the vector $(e(A_i))_i$, rather than a simple sum.
Formally, let
$$D^b_{k, p}(G) = \min_{A_1, \ldots, A_k} ||(e(A_i))_i||_p,$$
where the minimum is taken over all balanced $k$-partitions of $G$.
Note that $D^b_{k, 1} = D^b_k$, and
$D^b_{k, \infty}(G) = \min_{A_1, \ldots, A_k} \max_i e(A_i).$
The only parameter that has received some attention aside from $D^b_{k, 1}$ is $D^b_{k, \infty}$, as other cases seem to be outside the reach of the currently known methodology.
Intuitively, the problem of determining $D^b_{k, \infty}$ amounts to finding a balanced $k$-partition in which the densest part spans as few edges as possible.

\begin{theorem}[\cite{10problems}]
  If $G$ is a triangle-free graph on $n$ vertices, then 
  $$D^b_{2, \infty}(G) \leq \frac{n^2}{18} \qquad \text{and}\qquad  D^b_{3, \infty}(G) \leq \frac{n^2}{48} + o(n^2).$$
  Moreover, $D^b_{2, \infty}(K_{\frac{2n}{3},\frac{n}{3}}) = \frac{n^2}{18}$ and $D^b_{3, \infty}(K_{\frac{n}{2},\frac{n}{2}}) = \frac{n^2}{48}$.
\end{theorem}

Interestingly, it seems as if proving tight bounds on $D^b_{k, \infty}$ in the class of $K_4$-free graphs is noticeably more difficult, despite early optimistic conjectures.
\begin{conjecture}[\cite{10problems}]\label{first-conjecture}
  If $G$ is a $K_4$-free graph on $n$ vertices, then
  $$D^b_{2, \infty}(G) \leq \frac{n^2}{16}.$$
\end{conjecture}
\noindent\Cref{first-conjecture} was quickly disproven by Xu and Zhang in \cite{dense-halves-arxiv}.

\begin{theorem}[\cite{dense-halves-arxiv}]\label{first-lower-bound}
  There exists an infinite family of $K_4$-free graphs such that every graph $G$ in this family satisfies
  $$D^b_{2, \infty}(G) = \frac{37}{24^2}|G|^2 \approx 0.0642 \cdot |G|^2 >  0.0625 \cdot |G|^2 = \frac{|G|^2}{16}.$$
\end{theorem}

Another important observation made in \cite{dense-halves-arxiv} is that it seems difficult to obtain even relatively weak upper bounds on $D^b_{2, \infty}$ in the class of $K_4$-free graphs.

\begin{theorem}[\cite{dense-halves-arxiv}]\label{first-upper-bound}
  If $G$ is a $K_4$-free graph on $n>10^5$ vertices, then
  $$D^b_{2, \infty}(G) \leq 0.074 n^2.$$
\end{theorem}

In this work, we examine more closely the most relevant constructions for extremal values of $D^b_{2, \infty}$ in the class of $K_4$-free graphs and employ more advanced techniques to significantly narrow the gap between the bounds in \Cref{first-lower-bound} and \Cref{first-upper-bound}, from $0.0098$ to below $0.0017$.
First, we show a family of graphs which gives an even better lower bound.

\begin{theorem}\label{lower-bound-general}
  There exists an infinite family of $K_4$-free graphs, such that every graph $G$ in this family satisfies
  $$D^b_{2, \infty}(G) = \frac{425}{6498}|G|^2 \approx 0.0654 \cdot |G|^2.$$
\end{theorem}

Afterwards, using our improved understanding of potential extremal constructions, we leverage flag algebras to obtain a significantly better upper bound.
\begin{theorem}\label{upper-bound}
  If $G$ is a $K_4$-free graph on $n$ vertices, then
  $$D^b_{2, \infty}(G) \leq 0.0671n^2 + 5n^\frac{5}{3}.$$
\end{theorem}

%% file: sections/lower-bound.tex
Before we describe the family of graphs proving \Cref{lower-bound-general}, we recall some useful notions and a technical lemma from \cite{dense-halves-arxiv}.
An \emph{$n$-blowup} of a graph $G$ is a graph obtained by replacing every vertex in $G$ with $n$ copies.
Vertices in the blowup are connected if and only if their original counterparts in $G$ are connected.

\begin{lemma}[\cite{dense-halves-arxiv}]\label{edges-in-c5}
  If $F$ is an induced subgraph of order $pn + q$ of the $n$-blowup of $C_5$, where $p$ and $q$ are non-negative integers with $p < 5$ and $q < n$, then
  \[
    e(F) \geq f(n, p, q) := 
    \begin{cases}
      0, &\text{ if }\ p \leq 1, \\
      qn, &\text{ if }\ p = 2, \\
      n^2 + 2 qn, &\text{ if }\ p = 3, \\
      3n^2 + 2qn, &\text{ if }\ p = 4.
    \end{cases}
  \]
\end{lemma}

\begin{theorem}\label{lower-bound}
  Let $G$ be a graph on $57$ vertices consisting of a $5$-blowup of a $C_5$ and an independent set on the remaining $32$ vertices, with all edges between them.  
  If $G_{2n}$ is the $2n$-blowup of $G$, then
  $$D^b_{2, \infty}(G_{2n}) = 850n^2 = \frac{425}{6498}|G_{2n}|^2.$$
\end{theorem}

\begin{proof}
  Let $H$ be the unique copy of a $10n$-blowup of $C_5$ in $G_{2n}$ and let $(A, B)$ be any balanced bipartition of $V(G_{2n})$.
  If $H$ is fully contained in $A$ or $B$, then
  $$\max(e(A), e(B)) = 5 \cdot 10n \cdot 10n + 7n \cdot 50n = 850 n^2.$$
  This implies that $D^b_{2, \infty}(G_{2n}) \leq 850n^2$.
  Now, consider the case in which both sides of the bipartition contain a proper subset of $H$.
  Without loss of generality, we can assume that $A$ is the side containing the larger subset of $H$.
  This subset must be of size $10pn + q$ with $p = 2$ and $5n \leq q < 10n$ or with $3 \leq p < 5$ and $0 \leq q < 10n$.

  By \Cref{edges-in-c5}, the number of edges in $A$ is
  \begin{equation}\label{lower-bound-ineq}
    \begin{aligned}
      e(A) &= |A \cap V(H)|(|A| - |A \cap V(H)|) + e(A \cap V(H)) \\
           &\geq (10pn + q)(57n - 10pn - q) + f(10n, p, q) \\
           &= \begin{cases}
             740n^2 + q(27n - q), &\text{ if }\ p = 2,\ 5n \leq q < 10n, \\
             910n^2 + q(17n - q), &\text{ if }\ p = 3,\ 0 \leq q < 10n, \\
             980n^2 - q(3n + q), &\text{ if }\ p = 4,\ 0 \leq q < 10n.
           \end{cases}
    \end{aligned}
  \end{equation}
  Due to the basic properties of quadratic functions in a bounded range, we can conclude that $e(A)$ is always at least $850n^2$.
  This gives the desired lower bound on $D^b_{2, \infty}(G_{2n})$.

  We note that assigning $H$ to a single side of the partition is not the only way of attaining $D^b_{2, \infty}(G_{2n})$. 
  Denote the parts of $H$ as $A_1$, $A_2$, $A_3$, $A_4$, $A_5$, and consider a balanced bipartition $(A, B)$ for which $A_1, A_3 \subseteq A$, $A_2, A_4 \subseteq B$, and $A_5$ is split evenly.
  Notice that for such a partition, both $e(A) = e(B)$ and equality in \eqref{lower-bound-ineq} must hold.
  Therefore, we have $\max(e(A),e(B)) = 850n^2$.
\end{proof}

Our improvement on the lower bound over \cite{dense-halves-arxiv} stems from the fact that the equality $D^b_{2, \infty}(G_{2n}) = \max(e(A),e(B))$ is achieved for two different partitions rather than one.
Consider the class of graphs consisting of a balanced blowup of a $C_5$ and an independent set on the remaining vertices, with all edges between them.
Notice that, for any graph in this class, the values of $\max(e(A),e(B))$ for the two extremal partitions in \Cref{lower-bound} can be viewed as functions of the size of the independent set.
A simple analysis shows that these functions are of opposite monotonicity and the maximum of their pointwise minimum lies at the point where they are equal.
This point corresponds to $G_{2n}$.

Our other experiments on this problem strongly suggest that the class described in the previous paragraph is indeed the most difficult case for this problem.
Because of this, we believe that an extremal example is contained in it.
As a result, we consider $G_{2n}$ to be a natural candidate for an extremal example for the general problem of determining the dense half in $K_4$-free graphs.

\begin{conjecture}
  If $G$ is a $K_4$-free graph on $n$ vertices, then
  $$D^b_{2, \infty}(G) \leq \frac{425}{6498}n^2 + o(n^2).$$
\end{conjecture}

%% file: sections/flag-algebras.tex
To obtain a new upper bound on $D^b_{2, \infty}(G)$ in the class of $K_4$-free graphs, we leverage the powerful framework of flag algebras developed by Razborov (\cite{flag-algebras}).

Flag algebras have proven very useful in solving various problems related to finding partitions with few edges inside each part.
Some examples can be found in \cite{bipartite_large}, \cite{sparse-halves-best}, \cite{tri-sparse} and \cite{10problems}.
A typical flag algebraic application for these problems is built around picking a small substructure at random and defining the partition by splitting the vertices according to their connections to this structure.
This gives inequalities which can be directly plugged into the flag algebraic framework and can lead to strong bounds, usually (but not always) with the help of a semidefinite programming solver, such as \cite{csdp}.

It can be said that we are still in the early phases of exploring the applicability of expressing graph cuts with flag algebras, as ways to address key obstacles are still being found.
For instance, \cite{10problems} presents how this framework can be extended to analyze problems for balanced partitions, rather than arbitrary ones.
The authors of \cite{balanced-bipartite-distance-k4-free} show how flag algebras on colored graphs can be used to omit unfavorable structures.
In this work, we lift another limitation of this framework by showcasing how graph cuts and flag algebras can be used to argue about a non-algebraic function on the vector of numbers of edges inside the partition (though at the cost of introducing an asymptotically negligible error term).

A formal introduction to the flag-algebraic framework is beyond the scope of this paper.
For this purpose, we refer the interested reader to the original paper by Razborov \cite{flag-algebras}.

In the remainder of this section, we introduce the notation that we use when working with flags.
In this paper, we consider flags in multiple theories: 2-vertex-colored graphs, 2-vertex-colored partially oriented graphs, and 3-vertex-colored partially oriented graphs.
All flags we use are $K_4$-free, irrespective of coloring and orientation (that is, if we drop all colors and orientations, the graph should still be $K_4$-free).
In the case of partially oriented graphs, we allow orientation only between specific colors of vertices.
We leverage this fact to limit the number of flags during computer-assisted calculations.

We denote untyped flags by depicting the graph using round vertices.
For instance, this is the flag for the untyped triangle with two blue vertices and a single red vertex:
$$\Triangle{blue}{blue}{red}.$$
If we orient the edge between blue vertices, we get
$$\OTriangle{blue}{blue}{red}.$$

To denote a typed flag $(G, \theta)$ of type $\sigma$ from $\mathcal{A}^\sigma$, we depict $G$ and signify the injective mapping $\theta$ by writing the labels inside the vertices.
The vertices in $\operatorname{im}(\theta)$ are also square instead of round.
For example, this is a rooted polychromatic edge with a single unlabeled red vertex connected to the blue endpoint of the edge:
$$\EdgeL{blue}{red}{red}.$$

In our notation we also use some wildcards, which are syntactic sugar to shorten expressions where all possibilities of colors or orientations have to be considered.
If a flag contains white vertices, it represents the sum of all non-isomorphic flags obtained by coloring the white vertices in all possible ways.
For instance, in the theory of 2-colored graphs (where colors are red and blue), a flag with a single white vertex is shorthand for the sum of two flags
$$\Triangle{blue}{blue}{white} = \Triangle{blue}{blue}{red} + \Triangle{blue}{blue}{blue},$$
while in the theory of 3-colored flags (where colors are red, blue and cyan) it can be expanded to a sum of three summands
$$\Triangle{blue}{blue}{white} = \Triangle{blue}{blue}{red} + \Triangle{blue}{blue}{blue} + \Triangle{blue}{blue}{cyan}.$$
When more white vertices are involved, the length of the sum increases and isomorphisms start playing a role.
In the theory of 2-vertex-colored graphs, we have
$$\Triangle{white}{white}{white} = \Triangle{red}{red}{red} + \Triangle{blue}{red}{red} + \Triangle{blue}{blue}{red} + \Triangle{blue}{blue}{blue}.$$
When considering partially oriented graphs, we use an edge with a small dot as a template for any orientation, including the lack of it.
To be precise,
$$\AEdge{red}{red} = \Edge{red}{red} + \OEdge{red}{red}.$$
We also combine white vertices with orientation wildcards.
For instance, in the theory of 2-vertex-colored partially oriented graphs:
$$\EdgeBOA{red}{blue}{white} = \EdgeBOO{red}{blue}{blue} + \EdgeBOo{red}{blue}{blue} + \EdgeBON{red}{blue}{blue} + \EdgeBOO{red}{blue}{red} + \EdgeBOo{red}{blue}{red} + \EdgeBON{red}{blue}{red}.$$

We implicitly assume that if a flag contains objects that are not present in the currently considered theory, then its density is automatically zero and it can be safely ignored.

\subsection*{Linear product constraints}
In a typical flag algebraic proof, one shows that constraints of the form $f \geq 0$ hold for some $f \in \mathcal{A}^\sigma$ and some type $\sigma$.
We then apply the averaging operator and include a constraint of the form $\unlab{f} \geq 0$ in the instance of SDP that we need to solve.
A popular technique to increase the `strength' of the original inequality $f \geq 0$ is to include multiple constraints of the form $\unlab{f \cdot F} \geq 0$, where each $F$ is a single partially labeled graph whose labels constitute the type $\sigma$.

We employ this technique for some constraints that we use in our proofs, particularly of the form $\VertexM{red}{white} - \frac{1}{2} \geq 0$.
One can verify that a straightforward application of the averaging operator to such constraints (without the linear product extension) would yield relatively uninteresting inequalities.

%% file: sections/upper-bound.tex
The goal of this section is to prove \Cref{upper-bound}.
In the following theorems, we only consider graphs on an even number of vertices and prove a bound of $0.0671n^2 + 4n^\frac{5}{3}$.
The proof can be extended to graphs on an odd number of vertices by removing a single vertex and adding it back to either side of the partition.
This adds an extra error term of at most $n$, which is absorbed by the gap between the proven $4n^\frac{5}{3}$ and $5n^\frac{5}{3}$ in the exact formulation of \Cref{upper-bound}.

For a graph $G$ and two disjoint sets $A, B \subseteq V(G)$, by $E(A, B)$ we denote the set of edges with one endpoint in $A$ and the other in $B$, while $e(A, B)$ is the number of such edges.
Similarly, by $E(A)$ we denote the set of edges inside $A$ and by $e(A)$ the size of this set.
We use the convention that if $A = \emptyset$, then $e(A)/|A|^k = 0$ and $e(A,B)/|A|^k = 0$ for $k \in \{1, 2\}$. 
For brevity, a $2$-partition $V(G) = A \cup B$ is called a \emph{cut} and edges within $A$ or $B$ are called \emph{class-edges}.

We first show how graph cuts can be used to argue about asymptotic bounds on $D^b_{2, \infty}$.
This result is a known fact (see, for example, \cite{LLS13}), but we present its full proof to introduce proof techniques that we make heavy use of later.

\begin{theorem}\label{fully-random}
  Every graph $G$ on $n$ vertices satisfies $D^b_{2, \infty}(G) \leq \frac{e(G)}{4} + n^\frac{5}{3}$.
\end{theorem}
\begin{proof}
  Pick a balanced bipartition $(A, B)$ of $V(G)$ uniformly at random.
  For an edge $uv \in E(G)$, let $X_{uv}$ be the indicator random variable of the event $\{u, v\} \subseteq A$.
  It is clear that 
  $$\E[X_{uv}] = \frac{\binom{n - 2}{\frac{n}{2} - 2}}{\binom{n}{\frac{n}{2}}} = \frac{n - 2}{4(n - 1)},$$
  and, since $X_{uv}$ is an indicator variable,
  $$\Var[X_{uv}] = \E[X_{uv}] - \E[X_{uv}]^2 = \frac{(n-2)(3n-2)}{16(n - 1)^2}.$$
  Let $uv$ and $wx$ be two vertex-disjoint edges in $E(G)$.
  Then
  \[
    \E[X_{uv}X_{wx}] = \frac{\binom{n - 4}{\frac{n}{2} - 4}}{\binom{n}{\frac{n}{2}}} = \frac{(n - 4)(n - 6)}{16(n - 3)(n - 1)}.
  \]
  As a result,
  \[
    \Cov(X_{uv}, X_{wx}) = \frac{(n - 4)(n - 6)}{16(n - 3)(n - 1)} - \frac{(n - 2)^2}{16(n - 1)^2} = - \frac{2n^2-9n+6}{8(n-3)(n-1)^2} < 0.
  \]
  Now, consider a pair of edges $uv$ and $uw$, which share exactly one vertex.
  The formula for the expectation of the product of $X_{uv}$ and $X_{uw}$ changes to
  \[
    \E[X_{uv}X_{uw}] = \frac{\binom{n - 3}{\frac{n}{2} - 3}}{\binom{n}{\frac{n}{2}}} = \frac{n - 4}{8(n - 1)}.
  \]
  This gives us the following covariance
  \[
    \Cov(X_{uv}, X_{uw}) = \frac{n - 4}{8(n - 1)} - \frac{(n - 2)^2}{16(n - 1)^2} = \frac{n^2 - 6n + 4}{16(n - 1)^2}.
  \]
  Let $X$ be the random variable describing the number of class-edges inside of $A$.
  Clearly, $X = \sum_{uv \in E(G)} X_{uv}$, and so $\E[X] = \frac{n-2}{4(n-1)}e(G)$.    
  The variance of $X$ can be upper bounded as follows,
  \begin{equation}
    \begin{aligned}
      \Var(X) &= \sum_{uv \in E(G)} \Var(X_{uv}) + \sum_{\substack{uv,wx \in E(G) \\ \{u, v\} \cap \{w, x\} = \emptyset}} \Cov(X_{uv}, X_{wx}) +
      \sum_{\substack{uv, uw \in E(G) \\ v \neq w}} \Cov(X_{uv}, X_{uw}) \\
              &\leq \frac{(n-2)(3n-2)}{16(n - 1)^2} \cdot n^2 + \frac{n^2 - 6n + 4}{16(n - 1)^2} \cdot n^3 \\
              &\leq \frac{(n-1)(3n-3)}{16(n - 1)^2} \cdot n^2 + \frac{n^2 - 2n + 1}{16(n - 1)^2} \cdot n^3 \leq \frac{n^3}{4}.
    \end{aligned}
  \end{equation}
  Now, a straightforward application of Chebyshev's inequality gives
  \begin{equation}
    \P\left(X \geq \frac{e(G)}{4} + n^\frac{5}{3}\right) \leq \P\left(\left|X - \frac{n - 2}{4(n - 1)}e(G)\right| \geq n^{\frac{5}{3}}\right) \leq \frac{\frac{1}{4}n^3}{n^{\frac{10}{3}}} = \frac{1}{4n^{\frac{1}{3}}}.
  \end{equation}

  Observe that an identical inequality must hold for the random variable describing the number of class-edges inside $B$.
  Using a simple union bound, we can conclude that every graph must have a balanced bipartition $(A, B)$, such that $\max(e(A), e(B)) \leq \frac{e(G)}{4} + n^{\frac{5}{3}}$.
\end{proof}

\Cref{fully-random} immediately proves \Cref{upper-bound} for small graphs.

\begin{theorem}\label{small}
  Every graph $G$ on $n \leq 4096$ vertices satisfies $D_{2, \infty}^b(G) \leq 3n^{\frac{5}{3}}$.
\end{theorem}

\begin{proof}
By \Cref{fully-random}, 
  \[
    D_{2, \infty}^b(G) \leq \frac{e(G)}{4} + n^{\frac{5}{3}} 
    \leq \frac{n^2}{8} + n^{\frac{5}{3}}.
  \]
  Since for $n \leq 4096$, $n^2 \leq 16 n^{\frac{5}{3}}$, this implies that $D_{2, \infty}^b(G) \leq 3n^{\frac{5}{3}}$.
\end{proof}

We now generalize \Cref{fully-random} to more complex randomized balanced bipartitions.
\begin{theorem}\label{fix-expectation}
  Let $G$ be a graph on $n > 4096$ vertices and let $(A, B, C)$ be a (potentially unbalanced) tripartition of $V(G)$, such that $|A| \leq \frac{n}{2}$ and $|B| \leq \frac{n}{2}$. Then
\begin{align*}
    D^b_{2, \infty}(G) \leq \max &\left(e(A) + \frac{\frac{n}{2} - |A|}{|C|} \cdot e(A, C) + \left(\frac{\frac{n}{2} - |A|}{|C|}\right)^2 e(C),\right. \\
   &\ \ \left. e(B) + \frac{\frac{n}{2} - |B|}{|C|} \cdot e(B, C) + \left(\frac{\frac{n}{2} - |B|}{|C|}\right)^2 e(C)\right) + 2n^\frac{5}{3}.
\end{align*}
\end{theorem}

\begin{proof}
  If $|C| = 0$, then $(A, B)$ already constitutes a balanced bipartition proving the statement.

  First, we assume that $C$ has at least 4 vertices and that the sizes of $A$ and $B$ do not exceed $\frac{n}{2} - 2$.
  The remaining case will be addressed at the end of the proof.
  
  Let $(A', B')$ be a randomized balanced bipartition obtained by randomly extending $A$ and $B$ to sets of size $\frac{n}{2}$ with vertices from $C$.
  For any edge $uv \in E(G)$, let $X_{uv}$ be an indicator variable of the event $\{u, v\} \subseteq A'$.
  If $uv \in E(A, C)$, then $\E[X_{uv}] = \frac{\frac{n}{2} - |A|}{|C|}$.
  Whereas if $uv \in E(C)$, then
  $$\E[X_{uv}] = \frac{\binom{|C| - 2}{\frac{n}{2} - |A| - 2}}{\binom{|C|}{\frac{n}{2} - |A|}} = \frac{(\frac{n}{2} - |A| - 1) \cdot (\frac{n}{2} - |A|)}{(|C| - 1)|C|}.$$
  Now, similarly to the proof of \Cref{fully-random}, we want to analyze the covariance of pairs of edges that end up as class-edges in $A'$.

  \noindent\textbf{Option 1:} $uv, wx \in E(A, C)$, such that $u, w \in A$ and $v \neq x$.

  We have
  \[
    \E[X_{uv}X_{wx}] = \frac{\binom{|C| - 2}{\frac{n}{2} - |A| - 2}}{\binom{|C|}{\frac{n}{2} - |A|}} = \frac{(\frac{n}{2} - |A| - 1) \cdot (\frac{n}{2} - |A|)}{(|C| - 1) \cdot |C|},
  \]
  which allows us to prove the non-positivity of the covariance
\begin{align*}
	&\Cov(X_{uv}, X_{wx}) = \E[X_{uv}X_{wx}] - \E[X_{uv}]\E[X_{wx}] = \\
    =\ & \frac{(\frac{n}{2} - |A| - 1) \cdot (\frac{n}{2} - |A|)}{(|C| - 1) \cdot |C|} - \frac{(\frac{n}{2} - |A|)^2}{|C|^2} = \frac{\frac{n}{2} - |A|}{|C|} \cdot
    \left( \frac{\frac{n}{2} - |A| - 1}{|C| - 1} - \frac{\frac{n}{2} - |A|}{|C|} \right) \\
    =\ & \left(\frac{\frac{n}{2} - |A|}{|C|}\right) \cdot \frac{\frac{n}{2} - |A| - |C|}{(|C| - 1)|C|} \leq 0.
\end{align*}

  \noindent\textbf{Option 2:} $uv, uw \in E(A, C)$, where $u \in C$.

  It holds that
  \[
    \E[X_{uv}X_{uw}] = \frac{\binom{|C| - 1}{\frac{n}{2} - |A| - 1}}{\binom{|C|}{\frac{n}{2} - |A|}} = \frac{\frac{n}{2} - |A|}{|C|},
  \]
  so the covariance is at most
  \[
    \Cov(X_{uv}, X_{uw}) = \frac{\frac{n}{2} - |A|}{|C|} - \left(\frac{\frac{n}{2} - |A|}{|C|}\right)^2 \leq \frac{\frac{n}{2} - |A|}{|C|} \leq \frac{n}{|C|}.
  \]

  \noindent\textbf{Option 3:} $uv \in E(A, C)$, $wx \in E(C)$, such that $v \in C$, $v \neq w$ and $v \neq x$.

  In this case
  \[
    \E[X_{uv}X_{wx}] = \frac{\binom{|C| - 3}{\frac{n}{2} - |A| - 3}}{\binom{|C|}{\frac{n}{2} - |A|}} = \frac{(\frac{n}{2} - |A| - 2) \cdot (\frac{n}{2} - |A| - 1) \cdot (\frac{n}{2} - |A|)}{(|C| - 2) \cdot (|C| - 1) \cdot |C|},
  \]
  which gives the following covariance
\begin{align*}
   &\Cov(X_{uv}, X_{wx}) = \E[X_{uv}X_{wx}] - \E[X_{uv}]\E[X_{wx}] \\
   =\ & \frac{(\frac{n}{2} - |A| - 2) (\frac{n}{2} - |A| - 1) (\frac{n}{2} - |A|)}{(|C| - 2) (|C| - 1) |C|}
    - \frac{(\frac{n}{2} - |A| - 1) (\frac{n}{2} - |A|) ^ 2}{(|C| - 1)|C|^2} \\
   =\ & \left(\frac{(\frac{n}{2} - |A| - 1) (\frac{n}{2} - |A|)}{(|C| - 1)|C|}\right) \left( \frac{\frac{n}{2} - |A| - 2}{|C| - 2} - \frac{\frac{n}{2} - |A|}{|C|}\right) \\
   =\ & \left(\frac{(\frac{n}{2} - |A| - 1) (\frac{n}{2} - |A|)}{(|C| - 1)|C|}\right) \frac{n - 2|A| - 2|C|}{(|C| - 2)|C|} \leq 0.
\end{align*}

  \noindent\textbf{Option 4:} $uv \in E(A, C)$, $uw \in E(C)$.

  For such a pair we have
  \[
    \E[X_{uv}X_{uw}] = \frac{\binom{|C| - 2}{\frac{n}{2} - |A| - 2}}{\binom{|C|}{\frac{n}{2} - |A|}} = \frac{(\frac{n}{2} - |A| - 1) (\frac{n}{2} - |A|)}{(|C| - 1)|C|},
  \]
  which gives the following covariance
\begin{align*}
    & \Cov(X_{uv}, X_{uw}) = \E[X_{uv}X_{uw}] - \E[X_{uv}]\E[X_{uw}] \\    
    =\ & \frac{(\frac{n}{2} - |A| - 1) \cdot (\frac{n}{2} - |A|)}{(|C| - 1)|C|} - \frac{(\frac{n}{2} - |A| - 1) \cdot (\frac{n}{2} - |A|) ^ 2}{(|C| - 1)|C| ^ 2} \\
    =\ & \frac{(\frac{n}{2} - |A| - 1) \cdot (\frac{n}{2} - |A|)}{(|C| - 1)|C|} \cdot \frac{|A| + |C| - \frac{n}{2}}{|C|} \leq \frac{n^2}{|C|^2}.
\end{align*}

  \noindent\textbf{Option 5:} $uv \in E(C)$, $wx \in E(C)$, such that $uv$ and $wx$ are vertex-disjoint.

  An analysis of this pair gives
  \[
    \E[X_{uv}X_{wx}] = \frac{\binom{|C| - 4}{\frac{n}{2} - |A| - 4}}{\binom{|C|}{\frac{n}{2} - |A|}} = \frac{(\frac{n}{2} - |A| - 3) (\frac{n}{2} - |A| - 2) (\frac{n}{2} - |A| - 1) (\frac{n}{2} - |A|)}{(|C| - 3)(|C| - 2)(|C| - 1)|C|}.
  \]
  Fortunately, the covariance can be shown to be non-positive
\begin{align*}
    & \Cov(X_{uv}, X_{wx}) = \E[X_{uv}X_{wx}] - \E[X_{uv}]\E[X_{wx}] \\
    =\ & \frac{(\frac{n}{2} - |A| - 3) (\frac{n}{2} - |A| - 2) (\frac{n}{2} - |A| - 1) (\frac{n}{2} - |A|)}{(|C| - 3)(|C| - 2)(|C| - 1)|C|}
    - \frac{(\frac{n}{2} - |A| - 1)^2 (\frac{n}{2} - |A|)^2}{(|C| - 1)^2|C|^2} \\
    =\ & \left(\frac{(\frac{n}{2} - |A| - 1) (\frac{n}{2} - |A|)}{(|C| - 1)|C|}\right)
    \cdot \frac{\left(\frac{n}{2} - |A| - |C|\right) \cdot \left(\left(n - 2|A| - 3\right)(2|C| - 3) - 3\right)}{(|C| - 3)(|C| - 2)(|C| - 1)|C|} \leq 0.
\end{align*}

  \noindent\textbf{Option 6:} $uv \in E(C)$, $uw \in E(C)$.

  In this final case we have
  \[
    \E[X_{uv}X_{uw}] = \frac{\binom{|C| - 3}{\frac{n}{2} - |A| - 3}}{\binom{|C|}{\frac{n}{2} - |A|}} = \frac{(\frac{n}{2} - |A| - 2) (\frac{n}{2} - |A| - 1) (\frac{n}{2} - |A|)}{(|C| - 2) (|C| - 1) |C|},
  \]
  so the covariance is at most
\begin{align*}
   & \Cov(X_{uv}, X_{uw}) = \E[X_{uv}X_{uw}] - \E[X_{uv}]\E[X_{uw}] \\
   =\ & \frac{(\frac{n}{2} - |A| - 2) (\frac{n}{2} - |A| - 1) (\frac{n}{2} - |A|)}{(|C| - 2) (|C| - 1) |C|}
    - \frac{(\frac{n}{2} - |A| - 1)^2 (\frac{n}{2} - |A|)^2}{(|C| - 1)^2|C|^2} \\
   \leq\ & \frac{(\frac{n}{2} - |A| - 2) (\frac{n}{2} - |A| - 1) (\frac{n}{2} - |A|)}{(|C| - 2) (|C| - 1) |C|} \leq \frac{n^3}{|C|^3}.
\end{align*}

  We now set $X$ to be the random variable describing the number of class-edges inside~$A'$.
  Its variance equals
  \begin{equation}
    \begin{aligned}
    \Var(X) &= \sum_{uv \in E(G)} \Var(X_{uv}) + \sum_{\substack{uv, wx \in E(G)\\ uv \neq wx}} \Cov(X_{uv}, X_{wx}) \\
              &\leq \frac{\frac{n}{2} - |A|}{|C|} \cdot |A||C|
              + \frac{(\frac{n}{2} - |A| - 1) (\frac{n}{2} - |A|)}{(|C| - 1)|C|} \cdot |C|^2 \\
              &\overset{\text{option 2}}{+} \frac{n}{|C|} \cdot |A|^2|C| \overset{\text{option 4}}{+} 2\cdot\frac{n^2}{|C|^2} \cdot |A||C|^2 \overset{\text{option 6}}{+} \frac{n^3}{|C|^3} \cdot |C|^3 \\
              &\leq \frac{7n^2}{12} + \frac{9n^3}{4} \leq 3n^3,
    \end{aligned}
  \end{equation}
where the second sum runs over ordered pairs of distinct edges, so every unordered pair of edges in option 4 is counted twice.
  A simple calculation gives
\begin{align*}
    \E[X] &= e(A) + \frac{\frac{n}{2} - |A|}{|C|} \cdot e(A, C) + \left(\frac{(\frac{n}{2} - |A| - 1) (\frac{n}{2} - |A|)}{(|C| - 1)|C|}\right) e(C) \\
  &\leq e(A) + \frac{\frac{n}{2} - |A|}{|C|} \cdot e(A, C) + \left(\frac{\frac{n}{2} - |A|}{|C|}\right)^2 e(C),
\end{align*}
  and a straightforward application of Chebyshev's inequality (analogous to the proof of \Cref{fully-random}) shows that a deviation from this expectation larger than $n^{\frac{5}{3}}$ can only happen with probability at most $3n^{-\frac{1}{3}}$.
  An analogous bound holds for the number of class-edges inside $B'$.
  A union bound, together with the lower bound on $n$, implies that there always exists a balanced bipartition $(A', B')$ in which the number of class-edges in neither part exceeds its expectation by more than $n^{\frac{5}{3}}$.
  Combined with the upper bound on the expected value, this gives  
  \begin{equation}\label{expectation-first}
    \begin{aligned}
      \max(e(A'), e(B')) \leq \max&\left(e(A) + \frac{\frac{n}{2} - |A|}{|C|} \cdot e(A, C) + \left(\frac{\frac{n}{2} - |A|}{|C|}\right)^2 e(C),\right. \\
      &\ \ \left.e(B) + \frac{\frac{n}{2} - |B|}{|C|} \cdot e(B, C) + \left(\frac{\frac{n}{2} - |B|}{|C|}\right)^2 e(C)\right) + n^\frac{5}{3}.
    \end{aligned}
  \end{equation}
  
  Inequality \eqref{expectation-first} holds under the assumptions that $C$ has at least 4 vertices and the sizes of $A$ and $B$ do not exceed $\frac{n}{2} - 2$.
  These requirements can be satisfied by moving at most 2 vertices from $A$ to $C$ and at most 2 vertices from $B$ to $C$.
  Let $(\widehat A, \widehat B, \widehat C)$ be the resulting partition, for which inequality \eqref{expectation-first} holds.
  We have $e(\widehat A) \leq e(A)$.
  Since we move at most $4$ vertices to $C$, it holds that $e(\widehat A, \widehat C) \leq e(A, C) + 2n$ and $e(\widehat C) \leq e(C) + 4n$.
  We also observe that
  \[
    \frac{\frac{n}{2} - |\widehat A|}{|\widehat C|} - \frac{\frac{n}{2} - |A|}{|C|} \leq \frac{\frac{n}{2} - |\widehat A|}{|C|} - \frac{\frac{n}{2} - |A|}{|C|} = \frac{|A| - |\widehat A|}{|C|} \leq \frac{2}{|C|}.
  \]
  This allows us to establish
\begin{align*}
  \frac{\frac{n}{2} - |\widehat A|}{|\widehat C|} \cdot e(\widehat A, \widehat C) - \frac{\frac{n}{2} - |A|}{|C|} \cdot e(A, C) & \leq \left(\frac{\frac{n}{2} - |\widehat A|}{|\widehat C|} - \frac{\frac{n}{2} - |A|}{|C|}\right) \cdot e(A, C) + 2n \\
  & \leq \frac{2}{|C|} \cdot e(A, C) + 2n \leq 2|A| + 2n \leq 3n
\end{align*}
and
\begin{align*}
\left(\frac{\frac{n}{2} - |\widehat A|}{|\widehat C|}\right)^2 e(\widehat C) - \left(\frac{\frac{n}{2} - |A|}{|C|}\right)^2 e(C) &\leq \left(\left(\frac{\frac{n}{2} - |\widehat A|}{|\widehat C|}\right)^2 - \left(\frac{\frac{n}{2} - |A|}{|C|}\right)^2\right) e(C) + 4n \\
&\leq \frac{2}{|C|} \cdot 2 \cdot \frac{|C|^2}{2} + 4n \leq 6n.
\end{align*}
Therefore, together with inequality \eqref{expectation-first}, we obtain
\begin{align*}
    D^b_{2, \infty}(G) \leq \max &\left(e(A) + \frac{\frac{n}{2} - |A|}{|C|} \cdot e(A, C) + \left(\frac{\frac{n}{2} - |A|}{|C|}\right)^2 e(C),\right. \\
   &\ \ \left. e(B) + \frac{\frac{n}{2} - |B|}{|C|} \cdot e(B, C) + \left(\frac{\frac{n}{2} - |B|}{|C|}\right)^2 e(C)\right) + n^\frac{5}{3} + 9n.
\end{align*}
Since for $n>27$, we have $9n < n^\frac{5}{3}$, the statement follows.
\end{proof}

\begin{theorem}\label{fix-expectation-symmetric}
  Let $G$ be a graph on $n > 4096$ vertices and let $(A, B, C, D)$ be a (potentially unbalanced) partition of $V(G)$, such that $|C|$ is even, $|A| + |C| \leq \frac{n}{2}$ and $|B| + |C| \leq \frac{n}{2}$.
  Then, there exists a balanced bipartition $(C_1, C_2)$ of $C$, which satisfies all of
  \begin{itemize}
    \item $e(A, C_1) \leq \frac{1}{2}e(A, C) + \frac{1}{3}n^{\frac{5}{3}}$,
    \item $e(B, C_2) \leq \frac{1}{2}e(B, C) + \frac{1}{3}n^{\frac{5}{3}}$,
    \item $e(C_1) \leq \frac{1}{4} e(C) + \frac{1}{3}n^{\frac{5}{3}}$,
    \item $e(C_2) \leq \frac{1}{4} e(C) + \frac{1}{3}n^{\frac{5}{3}}$,
    \item $e(D, C_1) \leq \frac{e(C, D)}{2} + \frac{1}{3}n^{\frac{5}{3}}$, and
    \item $e(D, C_2) \leq \frac{e(C, D)}{2} + \frac{1}{3}n^{\frac{5}{3}}$
  \end{itemize}
  at once.
\end{theorem}

\begin{proof}
  First, let us bound the probability that the first condition is violated by a uniformly random balanced bipartition $(C_1, C_2)$ of $C$.
  Consider an edge $uv \in E(A, C)$, and let $X_{uv}$ be the indicator variable of the event $uv \in E(A, C_1)$.
  It follows that $\E[X_{uv}] = \frac{1}{2}$ and $\Var[X_{uv}] = \frac{1}{4}$.
  Now, we want to consider the covariance of the indicator variables of two edges from $E(A, C)$.
  There are two options to consider.

  \noindent\textbf{Option 1:} $uv, wx \in E(A, C)$, with $u, w \in A$ and $v \neq x$.

  We have
  \[
    \E[X_{uv}X_{wx}] = \frac{\binom{|C| - 2}{\frac{|C|}{2} - 2}}{\binom{|C|}{\frac{|C|}{2}}} = \frac{|C| - 2}{4|C| - 4},
  \]
  which quickly implies $\Cov(X_{uv}, X_{wx}) \leq 0$.

  \noindent\textbf{Option 2:} $uv, uw \in E(A, C)$, such that $u \in C$ and $v \neq w$.

  Since $X_{uv}X_{uw}$ is the indicator of the event $u \in C_1$, we have $\E[X_{uv}X_{uw}] = \frac{1}{2}$, which gives a covariance of exactly $\frac{1}{4}$.

  Similarly to the proof of \Cref{fully-random}, we now consider a random variable describing the number of edges between $A$ and $C_1$, $X = \sum_{uv \in E(A, C)} X_{uv}$.
  A crude upper bound using $|C| \leq \frac{n}{2}$ gives
  \begin{equation}
    \begin{aligned}
      \Var(X) &= \sum_{uv \in E(A, C)} \Var(X_{uv}) +
      \sum_{\substack{uv,uw \in E(A, C) \\ u \in C,\ v \neq w}} \Cov(X_{uv}, X_{uw}) \\
              &\leq \frac{e(A, C)}{4} + \frac{1}{4} \cdot |A|^2 \cdot |C| \leq \frac{n^2}{8} + \frac{n^3}{8} \leq \frac{n^3}{4}.
    \end{aligned}
  \end{equation}
  Using Chebyshev's inequality as in \Cref{fully-random}, we get an upper bound on the probability of $e(A, C_1)$ exceeding its expected value by at least $\frac{1}{3}n^{\frac{5}{3}}$ equal to
  \[
    \frac{\frac{n^3}{4}}{\left(\frac{1}{3} n^{\frac{5}{3}}\right)^2} = \frac{9}{4}n^{-\frac{1}{3}}.
  \]

  The same bound holds for $e(D, C_1)$.
  At the same time, a bound on the probability that $e(C_1)$ deviates too much can be obtained identically to \Cref{fully-random}, yielding a weaker constant of $\frac{9}{4}$ instead of $\frac{1}{4}$.
  To apply the union bound over all six conditions, we need
  \[
    6 \cdot \frac{9}{4} n^{-\frac{1}{3}} < 1,
  \]
  which holds for $n > 2460$.
\end{proof}

\subsection{Cut definitions}\label{cut-definitions}
  We now define the cuts necessary to prove \Cref{upper-bound}.
  We heavily rely on the idea from \cite{balanced-bipartite-distance-k4-free} of endowing the graph with additional structure in order to encode extra information.

  When defining inequalities that follow from the cuts, we use division, and denominators might be equal to 0 for some constructions.
  In the subsequent parts of the proof we use these inequalities after multiplying both sides by the common denominator.
  This means that they also (trivially) hold for these constructions.

  \begin{definition}
    Let $G$ be a $K_4$-free graph on $n$ vertices.
    By $L(G) \subseteq V(G)$ we denote the set of vertices with degree less than $\frac{n}{2}$.
    By $H(G) \subseteq V(G)$ we denote the set of vertices with degree at least $\frac{n}{2}$.
  \end{definition}
  We usually refer to vertices in $L(G)$ as low-degree vertices, and to vertices in $H(G)$ as high-degree vertices.
  Similarly, low-degree (respectively, high-degree) edges are edges with both endpoints in $L(G)$ (respectively, in $H(G)$).

  \begin{definition}[\cite{balanced-bipartite-distance-k4-free}]
    For a $K_4$-free graph $G$, we define its \emph{degree-colored variant} as a 2-vertex-colored graph, which is obtained by coloring the vertices in $L(G)$ blue and the vertices in $H(G)$ red.
  \end{definition}

  Later, we will extend this idea and give special meaning to edge orientation.

  \begin{definition}
    Let $G$ be a $K_4$-free graph on $n$ vertices and let $v \in V(G)$ be a low-degree vertex.
    The \emph{low-degree vertex cut} is the randomized balanced bipartition $(A, A^c)$ obtained by setting $N(v) \subseteq A$, and distributing the remaining vertices randomly.
  \end{definition}

  Note that the upper bound on degree is necessary here to ensure formal soundness of the partition.
  For a fixed vertex $v$, the expected numbers of edges inside $A$ and $A^c$ are at most
  \[
    e(N(v)) + \frac{\frac{n}{2} - |N(v)|}{n - |N(v)|} \cdot e(N(v), V(G) \setminus N(v)) + \left(\frac{\frac{n}{2} - |N(v)|}{n - |N(v)|}\right)^2 \cdot e(V(G) \setminus N(v)),
  \]
  and
  \[
    \left(\frac{\frac{n}{2}}{n - |N(v)|}\right)^2 \cdot e(V(G) \setminus N(v)).
  \]
  If we now consider the limits of $\E[e(A)]$ and $\E[e(A^c)]$ (normalized by $\binom{n}{2}$) in the infinite sequence of increasing blowups of the degree-colored variant of $G$, we get the following expressions in the flag algebra of type $\RootedVertex{blue}$:
  \[
    \VertexMM{blue}{white}{white} + \frac{\frac{1}{2} - \VertexM{blue}{white}}{\VertexN{blue}{white}} \VertexNM{blue}{white}{white}
    + \left(\frac{\frac{1}{2} - \VertexM{blue}{white}}{\VertexN{blue}{white}}\right)^2 \VertexNN{blue}{white}{white},
  \]
  and
  \[
    \left(\frac{\frac{1}{2}}{\VertexN{blue}{white}}\right)^2 \VertexNN{blue}{white}{white}.
  \]

  Our ultimate goal is to leverage this cut to show that no counterexample to \Cref{upper-bound} can exist.
  A typical approach, analogous to \cite{10problems} and \cite{balanced-bipartite-distance-k4-free}, would be to use the lower bound on the numbers of edges inside this cut to show a contradiction.
  However, in our problem, we are considering a non-algebraic function of the vector of numbers of class-edges, which ``leaks'' into the inequality:

  \[
    \max\left(\VertexMM{blue}{white}{white} + \frac{\frac{1}{2} - \VertexM{blue}{white}}{\VertexN{blue}{white}} \VertexNM{blue}{white}{white} +
    \left(\frac{\frac{1}{2} - \VertexM{blue}{white}}{\VertexN{blue}{white}}\right)^2 \VertexNN{blue}{white}{white},
    \left(\frac{\frac{1}{2}}{\VertexN{blue}{white}}\right)^2 \VertexNN{blue}{white}{white}\right)
    \geq c,
  \]
  where $c$ is a constant, and, as a result, it cannot be used in a typical flag-algebraic computer-assisted proof.
  Fortunately, this cut can be transformed into a set of linear inequalities in a simple extension of the degree-colored variant of $G$.

  \begin{definition}
    The \emph{low-degree vertex cut colored variant} of $G$ is a modification of the degree-colored variant of $G$, which is obtained as a result of the following process:
    if a low-degree vertex $v$ satisfies
    \begin{equation}\label{blue-vertex-raw}
      \begin{split}
        &e(N(v)) + \frac{\frac{n}{2} - |N(v)|}{n - |N(v)|} \cdot e(N(v), V(G) \setminus N(v)) + \left(\frac{\frac{n}{2} - |N(v)|}{n - |N(v)|}\right)^2 \cdot e(V(G) \setminus N(v)) \\
        \geq\ &
        \left(\frac{\frac{n}{2}}{n - |N(v)|}\right)^2 \cdot e(V(G) \setminus N(v)),
      \end{split}
    \end{equation}
    it remains blue,
    whereas if it satisfies
    \begin{equation}\label{cyan-vertex-raw}
      \begin{split}
        &e(N(v)) + \frac{\frac{n}{2} - |N(v)|}{n - |N(v)|} \cdot e(N(v), V(G) \setminus N(v)) + \left(\frac{\frac{n}{2} - |N(v)|}{n - |N(v)|}\right)^2 \cdot e(V(G) \setminus N(v)) \\
        <\ & \left(\frac{\frac{n}{2}}{n - |N(v)|}\right)^2 \cdot e(V(G) \setminus N(v)),
      \end{split}
    \end{equation}
    it is recolored to cyan.
  \end{definition}
  Now, because of \eqref{blue-vertex-raw}, we can deduce the following identity for any blue vertex in the limit of the infinite sequence of increasing blowups of the low-degree vertex cut colored variant of $G$:
  \begin{equation}\label{blue-vertex}
    \begin{split}
      &\max \left(\VertexMM{blue}{white}{white} + \frac{\frac{1}{2} - \VertexM{blue}{white}}{\VertexN{blue}{white}} \VertexNM{blue}{white}{white} +
      \left(\frac{\frac{1}{2} - \VertexM{blue}{white}}{\VertexN{blue}{white}}\right)^2 \VertexNN{blue}{white}{white},
      \left(\frac{\frac{1}{2}}{\VertexN{blue}{white}}\right)^2 \VertexNN{blue}{white}{white}\right) \\
      =\ &\VertexMM{blue}{white}{white} + \frac{\frac{1}{2} - \VertexM{blue}{white}}{\VertexN{blue}{white}} \VertexNM{blue}{white}{white} +
      \left(\frac{\frac{1}{2} - \VertexM{blue}{white}}{\VertexN{blue}{white}}\right)^2 \VertexNN{blue}{white}{white}.
    \end{split}
  \end{equation}
  Analogously, for cyan vertices we get (because of \eqref{cyan-vertex-raw}):
  \begin{equation}\label{cyan-vertex}
    \begin{split}
      &\max\left(\VertexMM{cyan}{white}{white} + \frac{\frac{1}{2} - \VertexM{cyan}{white}}{\VertexN{cyan}{white}} \VertexNM{cyan}{white}{white} +
      \left(\frac{\frac{1}{2} - \VertexM{cyan}{white}}{\VertexN{cyan}{white}}\right)^2 \VertexNN{cyan}{white}{white},
      \left(\frac{\frac{1}{2}}{\VertexN{cyan}{white}}\right)^2 \VertexNN{cyan}{white}{white}\right) \\
      =\ & \left(\frac{\frac{1}{2}}{\VertexN{cyan}{white}}\right)^2 \VertexNN{cyan}{white}{white}.
    \end{split}
  \end{equation}

  \begin{definition}
    Let $G$ be a $K_4$-free graph on $n$ vertices and let $v \in V(G)$ be a high-degree vertex.
    Let $C$ be the smaller of the sets $H(G)$ and $L(G)$.
    The \emph{high-degree vertex cut} is the following randomized balanced bipartition $(A, A^c)$:
    \begin{itemize}
      \item $(V(G) \setminus N(v)) \subseteq A$,
      \item $(N(v) \cap C) \subseteq A^c$,
      \item remaining vertices are distributed randomly.
    \end{itemize}
  \end{definition}

  For a fixed high-degree vertex $v$ in a $K_4$-free graph $G$, the expected numbers of class-edges inside $A$ and $A^c$ in the high-degree vertex cut based on $v$ are at most
  \begin{equation}\label{high-degree-vertex-cut-raw-left}
    \begin{split}
      &e(V(G) \setminus N(v)) + \frac{\frac{n}{2} - |V(G) \setminus N(v)|}{|N(v) \setminus C|} \cdot e(V(G) \setminus N(v), N(v) \setminus C) \\
      &+ \left(\frac{\frac{n}{2} - |V(G) \setminus N(v)|}{|N(v) \setminus C|}\right)^2 \cdot e(N(v) \setminus C),
    \end{split}
  \end{equation}
  and
  \begin{equation}\label{high-degree-vertex-cut-raw-right}
    \begin{split}
      &e(N(v) \cap C) + \frac{\frac{n}{2} - |N(v) \cap C|}{|N(v) \setminus C|}e(N(v) \cap C, N(v) \setminus C) \\
      &+ \left(\frac{\frac{n}{2} - |N(v) \cap C|}{|N(v) \setminus C|}\right)^2 e(N(v) \setminus C).
    \end{split}
  \end{equation}
  Unfortunately, in the language of flags, we cannot express the set $C$.
  Because of this, to find unique limits that are equivalent to \eqref{high-degree-vertex-cut-raw-left} and \eqref{high-degree-vertex-cut-raw-right}, we have to consider different possibilities for which one of $L(G)$ and $H(G)$ is larger.
  If high-degree vertices are more prevalent in $G$, then the limits of $\E[e(A)]$ and $\E[e(A^c)]$ are equivalent to the following expressions in the flag algebra of type $\RootedVertex{red}$:
  \begin{equation*}
    \VertexNN{red}{white}{white} + \frac{\frac{1}{2} - \VertexN{red}{white}}{\VertexM{red}{red}} \VertexNM{red}{red}{white}
    + \left(\frac{\frac{1}{2} - \VertexN{red}{white}}{\VertexM{red}{red}}\right)^2 \VertexMM{red}{red}{red},
  \end{equation*}
  and
  \begin{equation*}
    \VertexMM{red}{blue}{blue} + \frac{\frac{1}{2} - \VertexM{red}{blue}}{\VertexM{red}{red}} \VertexMM{red}{blue}{red}
    + \left(\frac{\frac{1}{2} - \VertexM{red}{blue}}{\VertexM{red}{red}}\right)^2 \VertexMM{red}{red}{red}.
  \end{equation*}
  If high-degree vertices constitute less than half of $G$'s vertices, then the limits change to
  \begin{equation*}
    \VertexNN{red}{white}{white} + \frac{\frac{1}{2} - \VertexN{red}{white}}{\VertexM{red}{blue}} \VertexNM{red}{blue}{white}
    + \left(\frac{\frac{1}{2} - \VertexN{red}{white}}{\VertexM{red}{blue}}\right)^2 \VertexMM{red}{blue}{blue},
  \end{equation*}
  and
  \begin{equation*}
    \VertexMM{red}{red}{red} + \frac{\frac{1}{2} - \VertexM{red}{red}}{\VertexM{red}{blue}} \VertexMM{red}{blue}{red}
    + \left(\frac{\frac{1}{2} - \VertexM{red}{red}}{\VertexM{red}{blue}}\right)^2 \VertexMM{red}{blue}{blue}.
  \end{equation*}

  \begin{definition}\label{high-degree-edge-cut-def}
    Let $G$ be a $K_4$-free graph on $n$ vertices and let $uv$ be a high-degree edge.
    The \emph{high-degree edge cut} is the following randomized balanced bipartition $(A, A^c)$:
    \begin{itemize}
      \item $N(u) \setminus N(v) \subseteq A$,
      \item $N(v) \setminus N(u) \subseteq A^c$,
      \item $V(G) \setminus (N(u) \cup N(v))$ is randomly split into two balanced parts, one of which is included in $A$,
      \item remaining vertices are distributed randomly.
    \end{itemize}
  \end{definition}

  Notice that the high-degree requirement is essential for the soundness of this definition~-- since $\deg(u) \geq \frac{n}{2}$, we know that $|N(v) \setminus N(u)| + \frac{1}{2} |V(G) \setminus (N(u) \cup N(v))| \leq \frac{n}{2}$, and symmetrically for $v$.
  Denote by $N(u, v) = N(u)\cap N(v)$ the common neighborhood of $u$ and $v$.
  For a fixed high-degree edge $uv$, $\E[e(A)]$ is at most
  \begin{equation*}
    \begin{split}
      &e(N(u) \setminus N(v)) + \frac{1}{2} e(N(u) \setminus N(v), V(G) \setminus (N(u) \cup N(v))) + \frac{1}{4} e(V(G) \setminus (N(u) \cup N(v))) \\
      &+ \frac{\frac{n}{2} - |N(u) \setminus N(v)| - \frac{1}{2} |V(G) \setminus (N(u) \cup N(v))|}{|N(u, v)|} \cdot e(N(u) \setminus N(v), N(u, v)) \\
      &+ \frac{\frac{n}{4} - \frac{1}{2}|N(u) \setminus N(v)| - \frac{1}{4}|V(G) \setminus (N(u) \cup N(v))|}{|N(u, v)|} \cdot e(V(G) \setminus (N(u) \cup N(v)), N(u, v)),
    \end{split}
  \end{equation*}
  while $\E[e(A^c)]$ is at most
  \begin{equation*}
    \begin{split}
      &e(N(v) \setminus N(u)) + \frac{1}{2} e(N(v) \setminus N(u), V(G) \setminus (N(v) \cup N(u))) + \frac{1}{4} e(V(G) \setminus (N(v) \cup N(u))) \\
      &+ \frac{\frac{n}{2} - |N(v) \setminus N(u)| - \frac{1}{2} |V(G) \setminus (N(v) \cup N(u))|}{|N(v, u)|} \cdot e(N(v) \setminus N(u), N(v, u)) \\
      &+ \frac{\frac{n}{4} - \frac{1}{2}|N(v) \setminus N(u)| - \frac{1}{4}|V(G) \setminus (N(v) \cup N(u))|}{|N(v, u)|} \cdot e(V(G) \setminus (N(v) \cup N(u)), N(v, u)).
    \end{split}
  \end{equation*}

  The limits of these expressions in the infinite sequence of increasing blowups of the degree-colored variant of $G$ are in the flag algebra of type $\RootedEdge{red}{red}$:
  \[
    \EdgeLaL{red}{red}{white}{white} + \frac{1}{2} \EdgeLN{red}{red}{white}{white} + \frac{1}{4} \EdgeNaN{red}{red}{white}{white} +
    \frac{\frac{1}{2} - \EdgeL{red}{red}{white} - \frac{1}{2}\EdgeN{red}{red}{white}}{\EdgeB{red}{red}{white}} \EdgeBL{red}{red}{white}{white} +
    \frac{\frac{1}{4} - \frac{1}{2}\EdgeL{red}{red}{white} - \frac{1}{4}\EdgeN{red}{red}{white}}{\EdgeB{red}{red}{white}} \EdgeBN{red}{red}{white}{white},
  \]
  and
  \[
    \EdgeRaR{red}{red}{white}{white} + \frac{1}{2} \EdgeRN{red}{red}{white}{white} + \frac{1}{4} \EdgeNaN{red}{red}{white}{white} +
    \frac{\frac{1}{2} - \EdgeR{red}{red}{white} - \frac{1}{2}\EdgeN{red}{red}{white}}{\EdgeB{red}{red}{white}} \EdgeBR{red}{red}{white}{white} +
    \frac{\frac{1}{4} - \frac{1}{2}\EdgeR{red}{red}{white} - \frac{1}{4}\EdgeN{red}{red}{white}}{\EdgeB{red}{red}{white}} \EdgeBN{red}{red}{white}{white}.
  \]

  This cut, analogously to the low-degree vertex cut, is difficult to apply directly in the form of linear flag constraints.
  Similarly, we address this by considering an extension of our original graph, where information is directly encoded in its components.

  \begin{definition}
    The \emph{high-degree edge cut oriented variant} of $G$ is a partially oriented, vertex-colored graph, which is obtained as a result of the following modification of the low-degree vertex cut colored variant of $G$:\\
    For every high-degree edge $uv$, consider the high-degree edge cut $(A, A^c)$ based on this edge.
    If the above upper bound on $\E[e(A)]$ is at least the above upper bound on $\E[e(A^c)]$, orient this edge from $v$ to $u$.
    Otherwise, orient it from $u$ to $v$.
  \end{definition}

  For any choice of labels, the following relation must hold in the limit of the infinite sequence of increasing blowups of a high-degree edge cut oriented variant of a $K_4$-free graph $G$:
  \begin{equation}
    \begin{split}\label{high-degree-edge-cut-basis}
      &\OEdgeRaR{red}{red}{white}{white} + \frac{1}{2} \OEdgeRN{red}{red}{white}{white} + \frac{1}{4} \OEdgeNaN{red}{red}{white}{white} +
      \frac{\frac{1}{2} - \OEdgeR{red}{red}{white} - \frac{1}{2}\OEdgeN{red}{red}{white}}{\OEdgeB{red}{red}{white}} \OEdgeBR{red}{red}{white}{white} +
      \frac{\frac{1}{4} - \frac{1}{2}\OEdgeR{red}{red}{white} - \frac{1}{4}\OEdgeN{red}{red}{white}}{\OEdgeB{red}{red}{white}} \OEdgeBN{red}{red}{white}{white} \\
      \geq\ &\OEdgeLaL{red}{red}{white}{white} + \frac{1}{2} \OEdgeLN{red}{red}{white}{white} + \frac{1}{4} \OEdgeNaN{red}{red}{white}{white} +
      \frac{\frac{1}{2} - \OEdgeL{red}{red}{white} - \frac{1}{2}\OEdgeN{red}{red}{white}}{\OEdgeB{red}{red}{white}} \OEdgeBL{red}{red}{white}{white} +
      \frac{\frac{1}{4} - \frac{1}{2}\OEdgeL{red}{red}{white} - \frac{1}{4}\OEdgeN{red}{red}{white}}{\OEdgeB{red}{red}{white}} \OEdgeBN{red}{red}{white}{white}.
    \end{split}
  \end{equation}

  \begin{definition}
    Let $G$ be a $K_4$-free graph on $n$ vertices and let $uv$ be a low-degree edge.
    The \emph{low-degree edge cut} is the randomized balanced bipartition $(A, A^c)$ obtained by setting $N(u) \setminus N(v) \subseteq A$, $N(v) \setminus N(u) \subseteq A^c$ and distributing the remaining vertices randomly.
  \end{definition}

  The deterministically assigned sets $N(u)\setminus N(v)$ and $N(v)\setminus N(u)$ have fewer than $\frac{n}{2}$ vertices, since $u$ and $v$ are low-degree vertices. 
  For a fixed low-degree edge $uv$, $\E[e(A)]$ is at most
  \begin{equation*}
    \begin{split}
      &e(N(u) \setminus N(v)) + \frac{\frac{n}{2} - |N(u) \setminus N(v)|}{|(V(G) \setminus (N(u) \cup N(v))) \cup N(u, v)|} \\
      &\cdot \left(e(N(u) \setminus N(v), V(G) \setminus (N(u) \cup N(v))) + e(N(u) \setminus N(v), N(u, v))\right) \\
      &+ \left(\frac{\frac{n}{2} - |N(u) \setminus N(v)|}{|(V(G) \setminus (N(u) \cup N(v))) \cup N(u, v)|}\right)^2 \\
      &\cdot \left(e(V(G) \setminus (N(u) \cup N(v))) + e(V(G) \setminus (N(u) \cup N(v)), N(u, v))\right),
    \end{split}
  \end{equation*}
  while $\E[e(A^c)]$ is at most
  \begin{equation*}
    \begin{split}
      &e(N(v) \setminus N(u)) + \frac{\frac{n}{2} - |N(v) \setminus N(u)|}{|(V(G) \setminus (N(v) \cup N(u))) \cup N(u, v)|} \\
      &\cdot \left(e(N(v) \setminus N(u), V(G) \setminus (N(v) \cup N(u))) + e(N(v) \setminus N(u), N(v, u))\right) \\
      &+ \left(\frac{\frac{n}{2} - |N(v) \setminus N(u)|}{|(V(G) \setminus (N(v) \cup N(u))) \cup N(u, v)|}\right)^2 \\
      &\cdot \left(e(V(G) \setminus (N(v) \cup N(u))) + e(V(G) \setminus (N(v) \cup N(u)), N(v, u))\right).
    \end{split}
  \end{equation*}
  This time, the limits are in the algebra of type $\RootedEdge{blue}{blue}$:
  \[
    \EdgeLaL{blue}{blue}{white}{white} + 
    \frac{\frac{1}{2} - \EdgeL{blue}{blue}{white}}{\EdgeN{blue}{blue}{white} + \EdgeB{blue}{blue}{white}} \left(\EdgeLN{blue}{blue}{white}{white} + \EdgeBL{blue}{blue}{white}{white}\right) +
    \left(\frac{\frac{1}{2} - \EdgeL{blue}{blue}{white}}{\EdgeN{blue}{blue}{white} + \EdgeB{blue}{blue}{white}}\right)^2 \left(\EdgeNaN{blue}{blue}{white}{white} + \EdgeBN{blue}{blue}{white}{white}\right),
  \]
  and
  \[
    \EdgeRaR{blue}{blue}{white}{white} + 
    \frac{\frac{1}{2} - \EdgeR{blue}{blue}{white}}{\EdgeN{blue}{blue}{white} + \EdgeB{blue}{blue}{white}} \left(\EdgeRN{blue}{blue}{white}{white} + \EdgeBR{blue}{blue}{white}{white}\right) +
    \left(\frac{\frac{1}{2} - \EdgeR{blue}{blue}{white}}{\EdgeN{blue}{blue}{white} + \EdgeB{blue}{blue}{white}}\right)^2 \left(\EdgeNaN{blue}{blue}{white}{white} + \EdgeBN{blue}{blue}{white}{white}\right).
  \]
  Similarly to the high-degree edge cut, we need to define an extension of $G$, in which the limits of the expected numbers of class-edges will be usable in the form of linear constraints.

  \begin{definition}
    The \emph{low-degree edge cut oriented variant} of $G$ is a partially oriented, vertex-colored graph, which is obtained as a result of the following modification of the degree-colored variant of $G$:\\
    For every low-degree edge $uv$, consider the low-degree edge cut $(A, A^c)$ based on this edge.
    If the above upper bound on $\E[e(A)]$ is at least the above upper bound on $\E[e(A^c)]$, orient this edge from $v$ to $u$.
    Otherwise, orient it from $u$ to $v$.
  \end{definition}

  In the limit of the infinite sequence of increasing blowups of a low-degree edge cut oriented variant of some $K_4$-free graph $G$, a relationship analogous to \eqref{high-degree-edge-cut-basis} must hold:
  \begin{equation}
    \begin{split}
      &\OEdgeRaR{blue}{blue}{white}{white} + 
      \frac{\frac{1}{2} - \OEdgeR{blue}{blue}{white}}{\OEdgeN{blue}{blue}{white} + \OEdgeB{blue}{blue}{white}} \left(\OEdgeRN{blue}{blue}{white}{white} + \OEdgeBR{blue}{blue}{white}{white}\right) +
      \left(\frac{\frac{1}{2} - \OEdgeR{blue}{blue}{white}}{\OEdgeN{blue}{blue}{white} + \OEdgeB{blue}{blue}{white}}\right)^2 \left(\OEdgeNaN{blue}{blue}{white}{white} + \OEdgeBN{blue}{blue}{white}{white}\right) \\
      \geq\ &
      \OEdgeLaL{blue}{blue}{white}{white} + 
      \frac{\frac{1}{2} - \OEdgeL{blue}{blue}{white}}{\OEdgeN{blue}{blue}{white} + \OEdgeB{blue}{blue}{white}} \left(\OEdgeLN{blue}{blue}{white}{white} + \OEdgeBL{blue}{blue}{white}{white}\right) +
      \left(\frac{\frac{1}{2} - \OEdgeL{blue}{blue}{white}}{\OEdgeN{blue}{blue}{white} + \OEdgeB{blue}{blue}{white}}\right)^2 \left(\OEdgeNaN{blue}{blue}{white}{white} + \OEdgeBN{blue}{blue}{white}{white}\right).
    \end{split}
  \end{equation}

  \begin{definition}
    Let $G$ be a $K_4$-free graph on $n$ vertices.
    Let $C$ be the smaller of the sets $H(G)$ and $L(G)$.
    The \emph{degree cut} is the randomized balanced bipartition $(A, A^c)$ obtained by setting $C \subseteq A$ and distributing the remaining vertices randomly.
  \end{definition}
  The expected numbers of class-edges inside $A$ and $A^c$ are at most
  \begin{equation*}
    e(C) + \frac{\frac{n}{2} - |C|}{|V(G) \setminus C|} \cdot e(C, V(G) \setminus C) + \left(\frac{\frac{n}{2} - |C|}{|V(G) \setminus C|}\right)^2 \cdot e(V(G) \setminus C),
  \end{equation*}
  and
  \begin{equation*}
    \left(\frac{\frac{n}{2}}{|V(G) \setminus C|}\right)^2 \cdot e(V(G) \setminus C).
  \end{equation*}

  If low-degree vertices are more prevalent in $G$, then the limits of $\E[e(A)]$ and $\E[e(A^c)]$ are the following expressions in the flag algebra of unrooted type:
  \begin{equation*}
    \Edge{red}{red} + \frac{\frac{1}{2} - \Vertex{red}}{\Vertex{blue}} \Edge{red}{blue} + \left(\frac{\frac{1}{2} - \Vertex{red}}{\Vertex{blue}}\right)^2 \Edge{blue}{blue},
    \text{ and }
    \left(\frac{\frac{1}{2}}{\Vertex{blue}}\right)^2 \Edge{blue}{blue}.
  \end{equation*}
  In the case where high-degree vertices are more common, we get the following expressions:
  \begin{equation*}
    \Edge{blue}{blue} + \frac{\frac{1}{2} - \Vertex{blue}}{\Vertex{red}} \Edge{blue}{red} + \left(\frac{\frac{1}{2} - \Vertex{blue}}{\Vertex{red}}\right)^2 \Edge{red}{red},
    \text{ and }
    \left(\frac{\frac{1}{2}}{\Vertex{red}}\right)^2 \Edge{red}{red}.
  \end{equation*}

  Notice that, contrary to the previous cuts, the expected number of class-edges for this cut is a global property of the graph, rather than a local one.
  As a result, we can work around the presence of the $\max$ function by separating graphs into classes based on where the maximum is achieved.
  For instance, if we limit ourselves to the class of graphs in which the following inequality holds:
  \begin{equation*}
    \begin{split}
      &e(C) + \frac{\frac{n}{2} - |C|}{|V(G) \setminus C|} \cdot e(C, V(G) \setminus C) 
      + \left(\frac{\frac{n}{2} - |C|}{|V(G) \setminus C|}\right)^2 \cdot e(V(G) \setminus C) \\
      \geq\ &
      \left(\frac{\frac{n}{2}}{|V(G) \setminus C|}\right)^2 \cdot e(V(G) \setminus C),
    \end{split}
  \end{equation*}
  we can freely use the degree cut, as in this class the following trivially holds:
  \begin{equation*}
    \begin{split}
      \max &\left(e(C) + \frac{\frac{n}{2} - |C|}{|V(G) \setminus C|} \cdot e(C, V(G) \setminus C) + \left(\frac{\frac{n}{2} - |C|}{|V(G) \setminus C|}\right)^2 \cdot e(V(G) \setminus C)\right., \\
           &\left.\left(\frac{\frac{n}{2}}{|V(G) \setminus C|}\right)^2 \cdot e(V(G) \setminus C)\right) \\
           =\ &
      e(C) + \frac{\frac{n}{2} - |C|}{|V(G) \setminus C|} \cdot e(C, V(G) \setminus C) + \left(\frac{\frac{n}{2} - |C|}{|V(G) \setminus C|}\right)^2 \cdot e(V(G) \setminus C).
    \end{split}
  \end{equation*}
  When considering the limits of graph sequences, a formally sound usage of the degree cut requires an analysis of four different cases as we have to account for the fact that there are two possibilities for what the set $C$ actually represents.

\subsection{Proof of the upper bound}
  We now combine all the cuts and related graph extensions defined above to prove a new upper bound on $D_{2, \infty}^b$ in the class of $K_4$-free graphs.

\begin{proof}[Proof of \Cref{upper-bound}]
  If $G$ has at most 4096 vertices, we are done by \Cref{small}, so we may assume that $n>4096$. As explained at the beginning of this section, we will prove the bound $D^b_{2, \infty}(G) < 0.0671 n^2 + 4n^\frac{5}{3}$ for even $n$, which implies the desired general bound $D^b_{2, \infty}(G) < 0.0671 n^2 + 5n^\frac{5}{3}$ for arbitrary $n$.

  As mentioned in the proof of \Cref{lower-bound}, we know two cuts which are tight for the conjectured extremal example.
  In the nomenclature used in this section, these are the high-degree edge cut and the degree cut.
  The remaining cuts are used to eliminate degenerate constructions which are so far from the extremal example that the optimal cuts are no longer relevant.

  We have observed that a formally sound application of the degree cut requires an analysis of four separate graph classes.
  In addition, we need to handle graphs with a small density of high-degree edges separately, as for them the high-degree edge cut gives no useful information.
  This, in total, gives 8 different cases.

  In each of these cases, we show that for at least one of the cuts defined in \Cref{cut-definitions} the maximum of expected values of the number of class-edges inside the two sides of the cut is below our upper bound.
  To do this, we consider the limit of the infinite sequence of increasing blowups of the low-degree (or high-degree) edge cut oriented variant of the graph,
  and leverage an SDP solver to show that a graphon for which every cut would ``exceed'' our upper bound cannot exist.

  We first outline a set of base inequalities that will be relevant in the computer-assisted flag algebraic part of every case.
  Our usage of colored vertices clearly implies
  \begin{equation}\label{degree-blue}
    \VertexM{blue}{white} \leq \frac{1}{2},
  \end{equation}
  \begin{equation}\label{degree-cyan}
    \VertexM{cyan}{white} \leq \frac{1}{2},
  \end{equation}
  \begin{equation}\label{degree-red}
    \VertexM{red}{white} \geq \frac{1}{2},
  \end{equation}
  and, as a consequence of the fully random partition,
  \begin{equation}\label{random-cut}
    \frac{1}{4}\Edge{white}{white} \geq 0.0671 \cdot 2.
  \end{equation}

  \textbf{Case 1:}
  In this case, we consider a $K_4$-free graph $G$ which satisfies the following set of inequalities:
  \begin{equation}\label{low-degree-prevalent-graph}
    \frac{|L(G)|}{n} \geq \frac{1}{2},
  \end{equation}
  \begin{equation}\label{high-degree-edges-present-graph}
    \frac{e(H(G))}{\binom{n}{2}} \geq 0.02,
  \end{equation}
  and
  \begin{equation}\label{degree-cut-case-1-graph}
    \begin{split}
      &e(H(G)) + \frac{\frac{n}{2} - |H(G)|}{|L(G)|} e(L(G), H(G)) + \left(\frac{\frac{n}{2} - |H(G)|}{|L(G)|}\right)^2 e(L(G)) \\
      \geq & \left(\frac{\frac{n}{2}}{|L(G)|}\right)^2 e(L(G)).
    \end{split}
  \end{equation}
  For contradiction, assume that $D^b_{2, \infty}(G) > 0.0671 n^2 + 4n^\frac{5}{3}$.

  Let $R$ be the high-degree edge cut oriented variant of $G$.
  We now consider the graphon $W$, which is the limit of the infinite sequence of increasing blowups of $R$.
  It is clear that, as a consequence of \eqref{low-degree-prevalent-graph}, \eqref{high-degree-edges-present-graph}, and \eqref{degree-cut-case-1-graph}, $W$ has to satisfy
  \begin{equation}\label{low-degree-prevalent}
    \Vertex{blue} + \Vertex{cyan} \geq \frac{1}{2},
  \end{equation}
  \begin{equation}\label{high-degree-edges-present}
    \Edge{red}{red} \geq 0.02,
  \end{equation}
  and
  \begin{equation}\label{degree-cut-case-1}
    \begin{split}
      &\Edge{red}{red} + \frac{\frac{1}{2} - \Vertex{red}}{\Vertex{blue} + \Vertex{cyan}} \left(\Edge{red}{blue} + \Edge{red}{cyan}\right) + \left(\frac{\frac{1}{2} - \Vertex{red}}{\Vertex{blue} + \Vertex{cyan}}\right)^2 \left(\Edge{blue}{blue} + \Edge{blue}{cyan} + \Edge{cyan}{cyan}\right) \\
      \geq\ &
      \left(\frac{\frac{1}{2}}{\Vertex{blue} + \Vertex{cyan}}\right)^2 \left(\Edge{blue}{blue} + \Edge{blue}{cyan} + \Edge{cyan}{cyan}\right).
    \end{split}
  \end{equation}
  For every blue vertex in $R$ it must hold that
  \begin{equation*}
    \begin{split}
      &e(N(v)) + \frac{\frac{n}{2} - |N(v)|}{n - |N(v)|} \cdot e(N(v), V(G) \setminus N(v)) + \left(\frac{\frac{n}{2} - |N(v)|}{n - |N(v)|}\right)^2 \cdot e(V(G) \setminus N(v)) \\
      \geq\ & 0.0671 n^2.
    \end{split}
  \end{equation*}
  After applying the standard normalization by $\binom{n}{2}$ and considering the limit, we get that the following must hold for $W$:
  \begin{equation}
    \VertexMM{blue}{white}{white} + \frac{\frac{1}{2} - \VertexM{blue}{white}}{\VertexN{blue}{white}} \VertexNM{blue}{white}{white} +
    \left(\frac{\frac{1}{2} - \VertexM{blue}{white}}{\VertexN{blue}{white}}\right)^2 \VertexNN{blue}{white}{white} \geq 0.0671 \cdot 2,
  \end{equation}
  which can be rewritten as
  \begin{equation}
    \left(\VertexMM{blue}{white}{white} - 0.1342\right) \cdot \VertexN{blue}{white}^2 + \left(\frac{1}{2} - \VertexM{blue}{white}\right) \cdot \VertexNM{blue}{white}{white} \cdot \VertexN{blue}{white} +
    \left(\frac{1}{2} - \VertexM{blue}{white}\right)^2 \VertexNN{blue}{white}{white}
    \geq 0.
  \end{equation}
  We can then apply the averaging operator and get a linear constraint of the form
  \begin{equation}\label{blue-vertex-sides}
    \unlab{
      \left(\VertexMM{blue}{white}{white} - 0.1342\right) \cdot \VertexN{blue}{white}^2 + \left(\frac{1}{2} - \VertexM{blue}{white}\right) \cdot \VertexNM{blue}{white}{white} \cdot \VertexN{blue}{white} +
      \left(\frac{1}{2} - \VertexM{blue}{white}\right)^2 \VertexNN{blue}{white}{white}
    }
    \geq 0.
  \end{equation}
  As a consequence of \eqref{blue-vertex}, $W$ must satisfy for any choice of the labeled blue vertex
  \begin{equation}
    \VertexMM{blue}{white}{white} + \frac{\frac{1}{2} - \VertexM{blue}{white}}{\VertexN{blue}{white}} \VertexNM{blue}{white}{white} +
    \left(\frac{\frac{1}{2} - \VertexM{blue}{white}}{\VertexN{blue}{white}}\right)^2 \VertexNN{blue}{white}{white} \geq
    \left(\frac{\frac{1}{2}}{\VertexN{blue}{white}}\right)^2 \VertexNN{blue}{white}{white},
  \end{equation}
  which can be multiplied and averaged to
  \begin{equation}\label{cut-blue-vertex}
    \unlab{
      \VertexMM{blue}{white}{white} \cdot \VertexN{blue}{white}^2 + \left(\frac{1}{2} - \VertexM{blue}{white}\right) \cdot \VertexNM{blue}{white}{white} \cdot \VertexN{blue}{white} -
      \left(\VertexM{blue}{white} - \VertexM{blue}{white}^2\right) \VertexNN{blue}{white}{white}
    }
    \geq 0.
  \end{equation}
  For cyan vertices, we get analogous constraints in the form of
  \begin{equation}\label{cyan-vertex-sides}
    \unlab{
      \frac{1}{4}\VertexNN{cyan}{white}{white} - 0.1342 \VertexN{cyan}{white}^2
    }
    \geq 0,
  \end{equation}
  and
  \begin{equation}\label{cyan-vertex-cut}
    \unlab{
      \left(\VertexM{cyan}{white} - \VertexM{cyan}{white}^2\right) \VertexNN{cyan}{white}{white} -
      \VertexMM{cyan}{white}{white} \cdot \VertexN{cyan}{white}^2 - \left(\frac{1}{2} - \VertexM{cyan}{white}\right) \cdot \VertexNM{cyan}{white}{white} \cdot \VertexN{cyan}{white}
    }
    \geq 0.
  \end{equation}
  By the assumptions of this case, we know $|L(G)| \geq |H(G)|$.
  For any high-degree vertex $v$ in $G$ we can also say
  \begin{equation}
    \begin{split}
      &e(V(G) \setminus N(v)) + \frac{\frac{n}{2} - |V(G) \setminus N(v)|}{|N(v) \setminus H(G)|} \cdot e(V(G) \setminus N(v), N(v) \setminus H(G)) \\
      &+ \left(\frac{\frac{n}{2} - |V(G) \setminus N(v)|}{|N(v) \setminus H(G)|}\right)^2 \cdot e(N(v) \setminus H(G)) \\
      \geq\ &
      e(N(v) \cap H(G)) + \frac{\frac{n}{2} - |N(v) \cap H(G)|}{|N(v) \setminus H(G)|} \cdot e(N(v) \cap H(G), N(v) \setminus H(G)) \\
      &+ \left(\frac{\frac{n}{2} - |N(v) \cap H(G)|}{|N(v) \setminus H(G)|}\right)^2 \cdot e(N(v) \setminus H(G)),
    \end{split}
  \end{equation}
  as the right side of the inequality represents the number of edges inside a subset of $N(v)$ of size $\frac{n}{2}$.
  By Mantel's theorem, this can never exceed $\frac{n^2}{16}$, so if the opposite inequality held for any high-degree $v$ in $G$, \Cref{fix-expectation} would immediately give a contradiction.
  As a consequence, we know that for any choice of the red label in $W$
  \begin{equation*}
    \begin{split}
      &\VertexNN{red}{white}{white} + \frac{\frac{1}{2} - \VertexN{red}{white}}{\VertexM{red}{blue} + \VertexM{red}{cyan}} \cdot \left(\VertexNM{red}{blue}{white} + \VertexNM{red}{cyan}{white}\right)
      + \left(\frac{\frac{1}{2} - \VertexN{red}{white}}{\VertexM{red}{blue} + \VertexM{red}{cyan}}\right)^2 \left(\VertexMM{red}{blue}{blue} + \VertexMM{red}{cyan}{blue} + \VertexMM{red}{cyan}{cyan}\right) \\
      \geq\ &
      \VertexMM{red}{red}{red} + \frac{\frac{1}{2} - \VertexM{red}{red}}{\VertexM{red}{blue} + \VertexM{red}{cyan}} \cdot \left(\VertexMM{red}{blue}{red} + \VertexMM{red}{cyan}{red}\right)
      + \left(\frac{\frac{1}{2} - \VertexM{red}{red}}{\VertexM{red}{blue} + \VertexM{red}{cyan}}\right)^2 \left(\VertexMM{red}{blue}{blue} + \VertexMM{red}{cyan}{blue} + \VertexMM{red}{cyan}{cyan}\right),
    \end{split}
  \end{equation*}
  which, after basic manipulation and averaging, implies
  \begin{equation}\label{high-degree-vertex-sides-case-1}
    \begin{split}
      &\left\llbracket\left(\VertexNN{red}{white}{white} - \VertexMM{red}{red}{red}\right) \cdot \left(\VertexM{red}{blue} + \VertexM{red}{cyan}\right)^2 \right.
      + \left(\frac{1}{2} - \VertexN{red}{white}\right) \cdot \left(\VertexM{red}{blue} + \VertexM{red}{cyan}\right) \cdot \left(\VertexNM{red}{blue}{white} + \VertexNM{red}{cyan}{white}\right) \\
      &- \left(\frac{1}{2} - \VertexM{red}{red}\right) \cdot \left(\VertexM{red}{blue} + \VertexM{red}{cyan}\right) \cdot \left(\VertexMM{red}{blue}{red} + \VertexMM{red}{cyan}{red}\right) \\
      &\left.+ \left[\left(\frac{1}{2} - \VertexN{red}{white}\right)^2 - \left(\frac{1}{2} - \VertexM{red}{red}\right)^2\right] \cdot \left(\VertexMM{red}{blue}{blue} + \VertexMM{red}{blue}{cyan} + \VertexMM{red}{cyan}{cyan}\right)\right\rrbracket \geq 0.
    \end{split}
  \end{equation}
  Since $D^b_{2, \infty}(G) > 0.0671 n^2 + 4n^\frac{5}{3}$, we also know that any high-degree vertex satisfies
  \begin{equation}
    \begin{split}
      &e(V(G) \setminus N(v)) + \frac{\frac{n}{2} - |V(G) \setminus N(v)|}{|N(v) \setminus H(G)|} \cdot e(V(G) \setminus N(v), N(v) \setminus H(G)) \\
      &+ \left(\frac{\frac{n}{2} - |V(G) \setminus N(v)|}{|N(v) \setminus H(G)|}\right)^2 \cdot e(N(v) \setminus H(G)) \geq 0.0671 n^2,
    \end{split}
  \end{equation}
  which means
  \begin{equation}\label{high-degree-vertex-cut-case-1}
    \begin{split}
      &\left\llbracket\left(\VertexNN{red}{white}{white} - 0.1342\right) \cdot \left(\VertexM{red}{blue} + \VertexM{red}{cyan}\right)^2\right.
      + \left(\frac{1}{2} - \VertexN{red}{white}\right) \cdot \left(\VertexM{red}{blue} + \VertexM{red}{cyan}\right) \cdot \left(\VertexNM{red}{blue}{white} + \VertexNM{red}{cyan}{white}\right) \\
      &\left.+ \left(\frac{1}{2} - \VertexN{red}{white}\right)^2 \cdot \left(\VertexMM{red}{blue}{blue} + \VertexMM{red}{blue}{cyan} + \VertexMM{red}{cyan}{cyan}\right)\right\rrbracket
      \geq 0.
    \end{split}
  \end{equation}
  We now consider the high-degree edge cut.
  We can reformulate, and average the assumption \eqref{high-degree-edge-cut-basis} into the form
  \begin{equation}\label{high-degree-edge-sides}
    \begin{split}
      &\left\llbracket\left(\OEdgeRaR{red}{red}{white}{white} - \OEdgeLaL{red}{red}{white}{white}\right) \cdot \OEdgeB{red}{red}{white} + \frac{1}{2} \left(\OEdgeRN{red}{red}{white}{white} - \OEdgeLN{red}{red}{white}{white}\right) \cdot \OEdgeB{red}{red}{white}\right. \\
      &+
      \left(\frac{1}{2} - \OEdgeR{red}{red}{white} - \frac{1}{2}\OEdgeN{red}{red}{white}\right) \cdot \OEdgeBR{red}{red}{white}{white} -
      \left(\frac{1}{2} - \OEdgeL{red}{red}{white} - \frac{1}{2}\OEdgeN{red}{red}{white}\right) \cdot \OEdgeBL{red}{red}{white}{white} \\
      &+
      \left.\left(\frac{1}{2}\OEdgeL{red}{red}{white} - \frac{1}{2}\OEdgeR{red}{red}{white}\right) \cdot \OEdgeBN{red}{red}{white}{white}\right\rrbracket
      \geq 0.
    \end{split}
  \end{equation}
  Similarly to the vertex cuts, we can see that $D^b_{2, \infty}(G) > 0.0671 n^2 + 4n^\frac{5}{3}$ implies
  \begin{equation}
    \begin{split}
      &\! \OEdgeRaR{red}{red}{white}{white} + \frac{1}{2} \OEdgeRN{red}{red}{white}{white} + \frac{1}{4} \OEdgeNaN{red}{red}{white}{white} +
      \frac{\frac{1}{2} - \OEdgeR{red}{red}{white} - \frac{1}{2}\OEdgeN{red}{red}{white}}{\OEdgeB{red}{red}{white}} \cdot \OEdgeBR{red}{red}{white}{white} +
      \frac{\frac{1}{4} - \frac{1}{2}\OEdgeR{red}{red}{white} - \frac{1}{4}\OEdgeN{red}{red}{white}}{\OEdgeB{red}{red}{white}} \cdot \OEdgeBN{red}{red}{white}{white} \\
      \geq\ & 0.0671 \cdot 2,
    \end{split}
  \end{equation}
  which, after averaging, implies
  \begin{equation}\label{high-degree-edge-cut}
    \begin{split}
      &\left\llbracket\left(\OEdgeRaR{red}{red}{white}{white} + \frac{1}{2} \OEdgeRN{red}{red}{white}{white} + \frac{1}{4} \OEdgeNaN{red}{red}{white}{white} - 0.1342\right) \cdot \OEdgeB{red}{red}{white}\right.
      + \left(\frac{1}{2} - \OEdgeR{red}{red}{white} - \frac{1}{2}\OEdgeN{red}{red}{white}\right) \cdot \OEdgeBR{red}{red}{white}{white} \\
      &\left.+ \left(\frac{1}{4} - \frac{1}{2}\OEdgeR{red}{red}{white} - \frac{1}{4}\OEdgeN{red}{red}{white}\right) \cdot \OEdgeBN{red}{red}{white}{white}\right\rrbracket \geq 0.
    \end{split}
  \end{equation}
  The degree-cut assumption \eqref{degree-cut-case-1} can be rewritten as
  \begin{equation}\label{degree-cut-sides-case-1}
    \begin{split}
      &\Edge{red}{red} \left(\Vertex{blue} + \Vertex{cyan}\right)^2 + \left(\frac{1}{2} - \Vertex{red}\right) \cdot \left(\Vertex{blue} + \Vertex{cyan}\right) \cdot \left(\Edge{red}{blue} + \Edge{red}{cyan}\right) \\
      &- \left(\Vertex{red} - \Vertex{red}^2 \right) \cdot \left(\Edge{blue}{blue} + \Edge{blue}{cyan} + \Edge{cyan}{cyan}\right) \geq 0.
    \end{split}
  \end{equation}
  It is also clear that the following must hold
  \begin{equation}
    \Edge{red}{red} + \frac{\frac{1}{2} - \Vertex{red}}{\Vertex{blue} + \Vertex{cyan}} \cdot \left(\Edge{red}{blue} + \Edge{red}{cyan}\right) + \left(\frac{\frac{1}{2} - \Vertex{red}}{\Vertex{blue} + \Vertex{cyan}}\right)^2 \cdot \left(\Edge{blue}{blue} + \Edge{blue}{cyan} + \Edge{cyan}{cyan}\right) \geq 0.0671 \cdot 2,
  \end{equation}
  which, in the linear form, is
  \begin{equation}\label{degree-cut-linear-case-1}
    \begin{split}
      &\left(\Edge{red}{red} - 0.1342\right) \cdot \left(\Vertex{blue} + \Vertex{cyan}\right)^2 + \left(\frac{1}{2} - \Vertex{red}\right) \cdot \left(\Vertex{blue} + \Vertex{cyan}\right) \cdot \left(\Edge{red}{blue} + \Edge{red}{cyan}\right) \\
                                                                                               &+ \left(\frac{1}{2} - \Vertex{red}\right)^2 \cdot \left(\Edge{blue}{blue} + \Edge{blue}{cyan} + \Edge{cyan}{cyan}\right) \geq 0.
    \end{split}
  \end{equation}
  The globality of this cut implies that no averaging is necessary.

  Calculations at \url{https://github.com/ignacy123/dense-half-k4-free-flags} show that no graphon can satisfy all of \eqref{degree-blue}, \eqref{degree-cyan}, \eqref{degree-red}, \eqref{random-cut},
  \eqref{low-degree-prevalent}, \eqref{high-degree-edges-present}, \eqref{blue-vertex-sides}, \eqref{cut-blue-vertex},
  \eqref{cyan-vertex-sides}, \eqref{cyan-vertex-cut}, \eqref{high-degree-vertex-sides-case-1}, \eqref{high-degree-vertex-cut-case-1},
  \eqref{high-degree-edge-sides}, \eqref{high-degree-edge-cut}, \eqref{degree-cut-sides-case-1} and \eqref{degree-cut-linear-case-1} at once.

  This clearly implies that the maximum of expected values of both sides of one of the cuts used in this case for $G$ is at most $0.0671 n^2$.
  Assume that the smallest cut is the high-degree edge cut for some edge $uv$.
  Let $(A, B, C, D)$ be the partition of $V(G)$ defined as $(N(u) \setminus N(v), N(v) \setminus N(u), V(G) \setminus (N(u) \cup N(v)), N(u) \cap N(v))$.
  If $|C|$ is odd, we choose any vertex $x$ from $C$ and apply \Cref{fix-expectation-symmetric} to $(A \cup \{x\}, B, C \setminus \{x\}, D)$ to get a split $(C_1, C_2)$ of $C \setminus \{x\}$.
  Afterwards, we obtain the final partition by applying \Cref{fix-expectation} to $(A \cup \{x\} \cup C_1, B \cup C_2, D)$, which can be done as $|A| + |C_1| + 1 \leq |A| + \frac{|C| + 1}{2} \leq |A| + |C| \leq \frac{n}{2}$.
  Our manipulations on $x$ introduce a slight error:
  \begin{itemize}
    \item The total number of edges introduced by moving $x$ to the left side of the partition can be trivially bounded by $\frac{n}{2}$,
    \item Edges with at least one endpoint in $N(u, v)$ are now more likely to end up on the right side of the partition than the formula following \Cref{high-degree-edge-cut-def} indicates.
      However, the increase of this probability is at most $\frac{1}{2|N(u, v)|}$.
      Crudely bounding the total number of edges with at least one endpoint in $N(u, v)$ by $n \cdot |N(u, v)|$, we can bound the total error by $\frac{n}{2}$.
  \end{itemize}
  Notice that $K_4$-freeness implies that there are no edges inside $N(u, v)$, which lets us conclude
  $D_{2, \infty}^b(G) \leq 0.0671 n^2 + \left(2 + \frac{1}{3} + \frac{1}{3} + \frac{1}{3}\right)n^\frac{5}{3} + n < 0.0671 n^2 + 4n^\frac{5}{3}$.
  If $|C|$ is even, an analogous but more straightforward application of \Cref{fix-expectation-symmetric} and \Cref{fix-expectation} gives the same bound but with a smaller error term.
  For other cuts, we only need \Cref{fix-expectation}, which gives $D_{2, \infty}^b(G) \leq 0.0671 n^2 + 2n^\frac{5}{3}$.
  In either case, we obtain a contradiction, giving a bound of $0.0671n^2 + 4n^\frac{5}{3}$ for even-sized graphs.
  The extension to odd-sized graphs described at the start of this section increases the error term to $5n^\frac{5}{3}$.

  \textbf{Case 2:}
  This case is almost identical to Case 1, but we consider the opposite relationship for the degree cut.
  To be precise, we consider $K_4$-free graphs that satisfy the following inequalities
  \begin{equation*}
    \frac{|L(G)|}{n} \geq \frac{1}{2},
  \end{equation*}
  \begin{equation*}
    \frac{e(H(G))}{\binom{n}{2}} \geq 0.02,
  \end{equation*}
  and
  \begin{equation}\label{degree-cut-case-2-graph}
    \begin{split}
      &\left(\frac{\frac{n}{2}}{|L(G)|}\right)^2 \cdot e(L(G)) \\
      \geq\ &
      e(H(G)) + \frac{\frac{n}{2} - |H(G)|}{|L(G)|} \cdot e(L(G), H(G)) + \left(\frac{\frac{n}{2} - |H(G)|}{|L(G)|}\right)^2 \cdot e(L(G)).
    \end{split}
  \end{equation}
  These conditions have the same implications for $W$ as in Case 1, aside from \eqref{degree-cut-case-1}, which changes to
  \begin{equation}\label{degree-cut-case-2}
    \begin{split}
      &\left(\frac{\frac{1}{2}}{\Vertex{blue} + \Vertex{cyan}}\right)^2 \left(\Edge{blue}{blue} + \Edge{blue}{cyan} + \Edge{cyan}{cyan}\right) \\
      \geq\ &
      \Edge{red}{red} + \frac{\frac{1}{2} - \Vertex{red}}{\Vertex{blue} + \Vertex{cyan}} \left(\Edge{red}{blue} + \Edge{red}{cyan}\right) + \left(\frac{\frac{1}{2} - \Vertex{red}}{\Vertex{blue} + \Vertex{cyan}}\right)^2 \left(\Edge{blue}{blue} + \Edge{blue}{cyan} + \Edge{cyan}{cyan}\right).
    \end{split}
  \end{equation}
  The linear form of \eqref{degree-cut-case-2} is
  \begin{equation}\label{degree-cut-sides-case-2}
    \begin{split}
      &\left(\Vertex{red} - \Vertex{red}^2 \right) \cdot \left(\Edge{blue}{blue} + \Edge{blue}{cyan} + \Edge{cyan}{cyan}\right) - \Edge{red}{red} \left(\Vertex{blue} + \Vertex{cyan}\right)^2 \\
      &- \left(\frac{1}{2} - \Vertex{red}\right) \cdot \left(\Vertex{blue} + \Vertex{cyan}\right) \cdot \left(\Edge{red}{blue} + \Edge{red}{cyan}\right) \geq 0.
    \end{split}
  \end{equation}
  By considering the denser side of this cut, we also get
  \begin{equation}
    \left(\frac{\frac{1}{2}}{\Vertex{blue} + \Vertex{cyan}}\right)^2 \left(\Edge{blue}{blue} + \Edge{blue}{cyan} + \Edge{cyan}{cyan}\right) \geq 0.0671 \cdot 2,
  \end{equation}
  which is the same as
  \begin{equation}\label{degree-cut-linear-case-2}
    \frac{1}{4} \left(\Edge{blue}{blue} + \Edge{blue}{cyan} + \Edge{cyan}{cyan}\right) - 0.1342\left(\Vertex{blue} + \Vertex{cyan}\right)^2  \geq 0.
  \end{equation}

  Calculations at \url{https://github.com/ignacy123/dense-half-k4-free-flags} show that no graphon can satisfy all of \eqref{degree-blue}, \eqref{degree-cyan}, \eqref{degree-red}, \eqref{random-cut},
  \eqref{low-degree-prevalent}, \eqref{high-degree-edges-present}, \eqref{blue-vertex-sides}, \eqref{cut-blue-vertex},
  \eqref{cyan-vertex-sides}, \eqref{cyan-vertex-cut}, \eqref{high-degree-vertex-sides-case-1}, \eqref{high-degree-vertex-cut-case-1},
  \eqref{high-degree-edge-sides}, \eqref{high-degree-edge-cut}, \eqref{degree-cut-sides-case-2} and \eqref{degree-cut-linear-case-2} at once.
  This implies a contradiction similarly to Case~1.

  \textbf{Case 3:}
  This case is also similar to Case 1, but it differs in which of $L(G)$ and $H(G)$ is larger.
  The exact set of inequalities which $G$ has to satisfy in this case is
  \begin{equation}\label{high-degree-prevalent-graph}
    \frac{|H(G)|}{n} \geq \frac{1}{2},
  \end{equation}
  \begin{equation}
    \frac{e(H(G))}{\binom{n}{2}} \geq 0.02,
  \end{equation}
  and
  \begin{equation}
    \begin{split}
      &e(L(G)) + \frac{\frac{n}{2} - |L(G)|}{|H(G)|} \cdot e(H(G), L(G)) + \left(\frac{\frac{n}{2} - |L(G)|}{|H(G)|}\right)^2 \cdot e(H(G)) \\
      \geq\ & \left(\frac{\frac{n}{2}}{|H(G)|}\right)^2 \cdot e(H(G)).
    \end{split}
  \end{equation}
  Clearly, $W$ must satisfy
  \begin{equation}\label{high-degree-prevalent}
    \Vertex{red} \geq \frac{1}{2},
  \end{equation}
  and
  \begin{equation}\label{degree-cut-raw-case-3}
    \begin{split}
      &\Edge{blue}{blue} + \Edge{blue}{cyan} + \Edge{cyan}{cyan} + \frac{\frac{1}{2} - \Vertex{blue} - \Vertex{cyan}}{\Vertex{red}} \left(\Edge{blue}{red} + \Edge{cyan}{red}\right)
      + \left(\frac{\frac{1}{2} - \Vertex{blue} - \Vertex{cyan}}{\Vertex{red}}\right)^2 \Edge{red}{red} \\
      \geq\ & \left(\frac{\frac{1}{2}}{\Vertex{red}}\right)^2 \Edge{red}{red}.
    \end{split}
  \end{equation}
  A crucial consequence of \eqref{high-degree-prevalent} is that \eqref{high-degree-vertex-sides-case-1} changes to
  \begin{equation}\label{high-degree-vertex-sides-case-3}
    \begin{split}
      &\left\llbracket\left(\VertexNN{red}{white}{white} - \VertexMM{red}{blue}{blue} - \VertexMM{red}{blue}{cyan} - \VertexMM{red}{cyan}{cyan}\right) \cdot \VertexM{red}{red}^2\right.
      + \left(\frac{1}{2} - \VertexN{red}{white}\right) \cdot \VertexM{red}{red} \cdot \VertexNM{red}{red}{white} \\
      &- \left(\frac{1}{2} - \VertexM{red}{blue} - \VertexM{red}{cyan}\right) \cdot \VertexM{red}{red} \cdot \left(\VertexMM{red}{red}{blue} + \VertexMM{red}{red}{cyan}\right) \\
      &\left.+ \left[\left(\frac{1}{2} - \VertexN{red}{white}\right)^2 - \left(\frac{1}{2} - \VertexM{red}{blue} - \VertexM{red}{cyan}\right)^2\right] \cdot \VertexMM{red}{red}{red}\right\rrbracket
        \geq 0,
    \end{split}
  \end{equation}
  and \eqref{high-degree-vertex-cut-case-1} changes to
  \begin{equation}\label{high-degree-vertex-cut-case-3}
    \begin{split}
    \left\llbracket\left(\VertexNN{red}{white}{white} - 0.1342\right) \cdot \VertexM{red}{red}^2 + \left(\frac{1}{2} - \VertexN{red}{white}\right) \cdot \VertexM{red}{red} \cdot \VertexNM{red}{red}{white}\right.
      \left.+ \left(\frac{1}{2} - \VertexN{red}{white}\right)^2 \VertexMM{red}{red}{red}\right\rrbracket
       \geq 0.
    \end{split}
  \end{equation}
  In addition, we have to rewrite \eqref{degree-cut-sides-case-1} as
  \begin{equation}\label{degree-cut-sides-case-3}
    \begin{split}
      &\left(\Edge{blue}{blue} + \Edge{blue}{cyan} + \Edge{cyan}{cyan}\right) \Vertex{red}^2 + \left(\frac{1}{2} - \Vertex{blue} - \Vertex{cyan}\right) \cdot \Vertex{red} \cdot \left(\Edge{blue}{red} + \Edge{cyan}{red}\right) \\
      &- \left(\Vertex{blue} + \Vertex{cyan} - (\Vertex{blue} + \Vertex{cyan})^2\right) \cdot \Edge{red}{red} \geq 0,
    \end{split}
  \end{equation}
  and \eqref{degree-cut-linear-case-1} as
  \begin{equation}\label{degree-cut-linear-case-3}
    \begin{split}
      &\left(\Edge{blue}{blue} + \Edge{blue}{cyan} + \Edge{cyan}{cyan} - 0.1342\right) \cdot \Vertex{red}^2 + \left(\frac{1}{2} - \Vertex{blue} - \Vertex{cyan}\right) \cdot \Vertex{red} \cdot \left(\Edge{blue}{red} + \Edge{cyan}{red}\right) \\
      &+ \left(\frac{1}{2} - \Vertex{blue} - \Vertex{cyan}\right)^2 \cdot \Edge{red}{red} \geq 0.
    \end{split}
  \end{equation}

  Calculations at \url{https://github.com/ignacy123/dense-half-k4-free-flags} show that no graphon can satisfy all of \eqref{degree-blue}, \eqref{degree-cyan}, \eqref{degree-red}, \eqref{random-cut},
  \eqref{high-degree-prevalent}, \eqref{high-degree-edges-present}, \eqref{blue-vertex-sides}, \eqref{cut-blue-vertex},
  \eqref{cyan-vertex-sides}, \eqref{cyan-vertex-cut}, \eqref{high-degree-vertex-sides-case-3}, \eqref{high-degree-vertex-cut-case-3},
  \eqref{high-degree-edge-sides}, \eqref{high-degree-edge-cut}, \eqref{degree-cut-sides-case-3} and \eqref{degree-cut-linear-case-3} at once.
  This implies a contradiction similarly to Case~1.

  \textbf{Case 4:}
  Case 4 is identical to Case 3 aside from the relationship between sides of the degree cut.
  Namely, we consider the following set of inequalities:

  \begin{equation*}
    \frac{|H(G)|}{n} \geq \frac{1}{2},
  \end{equation*}
  \begin{equation*}
    \frac{e(H(G))}{\binom{n}{2}} \geq 0.02,
  \end{equation*}
  and
  \begin{equation}
    \begin{split}
      &\left(\frac{\frac{n}{2}}{|H(G)|}\right)^2 \cdot e(H(G)) \\
      &\geq
      e(L(G)) + \frac{\frac{n}{2} - |L(G)|}{|H(G)|} \cdot e(H(G), L(G)) + \left(\frac{\frac{n}{2} - |L(G)|}{|H(G)|}\right)^2 \cdot e(H(G)).
    \end{split}
  \end{equation}
  Instead of \eqref{degree-cut-raw-case-3}, $W$ must now satisfy
  \begin{equation*}
    \begin{split}
      \left(\frac{\frac{1}{2}}{\Vertex{red}}\right)^2 \Edge{red}{red}
      \geq
      \Edge{blue}{blue}\! +\! \Edge{blue}{cyan}\! + \!\Edge{cyan}{cyan}\! + \frac{\frac{1}{2} - \Vertex{blue} - \Vertex{cyan}}{\Vertex{red}} \left(\Edge{blue}{red}\! + \!\Edge{cyan}{red}\right)
      + \left(\frac{\frac{1}{2} - \Vertex{blue} - \Vertex{cyan}}{\Vertex{red}}\right)^2 \Edge{red}{red},
    \end{split}
  \end{equation*}
  so \eqref{degree-cut-sides-case-3} changes to
  \begin{equation}\label{degree-cut-sides-case-4}
    \begin{split}
      &\left(\Vertex{blue} + \Vertex{cyan} - (\Vertex{blue} + \Vertex{cyan})^2 \right) \cdot \Edge{red}{red}
      - \left(\Edge{blue}{blue} + \Edge{blue}{cyan} + \Edge{cyan}{cyan}\right) \cdot \Vertex{red}^2 \\
      &- \left(\frac{1}{2} - \Vertex{blue} - \Vertex{cyan}\right) \cdot \Vertex{red} \cdot \left(\Edge{blue}{red} + \Edge{cyan}{red}\right) \geq 0,
    \end{split}
  \end{equation}
  and \eqref{degree-cut-linear-case-3} must change to
  \begin{equation}\label{degree-cut-linear-case-4}
    \begin{split}
      \frac{1}{4} \Edge{red}{red} - 0.1342 \Vertex{red}^2
      \geq 0.
    \end{split}
  \end{equation}

  Calculations at \url{https://github.com/ignacy123/dense-half-k4-free-flags} show that no graphon can satisfy all of \eqref{degree-blue}, \eqref{degree-cyan}, \eqref{degree-red}, \eqref{random-cut},
  \eqref{high-degree-prevalent}, \eqref{high-degree-edges-present}, \eqref{blue-vertex-sides}, \eqref{cut-blue-vertex},
  \eqref{cyan-vertex-sides}, \eqref{cyan-vertex-cut}, \eqref{high-degree-vertex-sides-case-3}, \eqref{high-degree-vertex-cut-case-3},
  \eqref{high-degree-edge-sides}, \eqref{high-degree-edge-cut}, \eqref{degree-cut-sides-case-4} and \eqref{degree-cut-linear-case-4} at once.
  This implies a contradiction similarly to Case~1.

  We now consider graphs whose density of high-degree edges is small.
  Aside from this change, the structure of Cases 5 to 8 mirrors the first four cases.

  \textbf{Case 5:} 
  We consider a $K_4$-free graph $G$ that satisfies
  \begin{equation*}
    \frac{|L(G)|}{n} \geq \frac{1}{2},
  \end{equation*}
  \begin{equation}\label{high-degree-edges-absent-graph}
    \frac{e(H(G))}{\binom{n}{2}} \leq 0.02,
  \end{equation}
  and
  \begin{equation*}
    \begin{split}
      &e(H(G)) + \frac{\frac{n}{2} - |H(G)|}{|L(G)|} \cdot e(L(G), H(G)) + \left(\frac{\frac{n}{2} - |H(G)|}{|L(G)|}\right)^2 \cdot e(L(G)) \\
      \geq\ & \left(\frac{\frac{n}{2}}{|L(G)|}\right)^2 \cdot e(L(G)).
    \end{split}
  \end{equation*}

  Now, by $R$ we denote the low-degree (rather than high-degree) edge cut oriented variant of $G$.
  Analogously, $W$ is the limit of the infinite sequence of increasing blowups of $R$.
  The most important implication of this change is that the low-degree vertex cut and the high-degree edge cut are no longer relevant.
  However, we can now use the low-degree edge cut.
  The analogue of \eqref{high-degree-edge-sides} for the low-degree edge cut is
  \begin{equation}\label{low-degree-edge-sides}
    \begin{split}
      &\left\llbracket\left(\OEdgeRaR{blue}{blue}{white}{white} - \OEdgeLaL{blue}{blue}{white}{white}\right) \cdot \left(\OEdgeN{blue}{blue}{white} + \OEdgeB{blue}{blue}{white}\right)^2\right. \\
      &+ \left(\frac{1}{2} - \OEdgeR{blue}{blue}{white}\right) \cdot \left(\OEdgeN{blue}{blue}{white} + \OEdgeB{blue}{blue}{white}\right) \cdot \left(\OEdgeRN{blue}{blue}{white}{white} + \OEdgeBR{blue}{blue}{white}{white}\right) \\
      &- \left(\frac{1}{2} - \OEdgeL{blue}{blue}{white}\right) \cdot \left(\OEdgeN{blue}{blue}{white} + \OEdgeB{blue}{blue}{white}\right) \cdot \left(\OEdgeLN{blue}{blue}{white}{white} + \OEdgeBL{blue}{blue}{white}{white}\right) \\
      &+ \left[\left(\frac{1}{2} - \OEdgeR{blue}{blue}{white}\right)^2 - \left(\frac{1}{2} - \OEdgeL{blue}{blue}{white}\right)^2\right] 
      \cdot \left.\left(\OEdgeNaN{blue}{blue}{white}{white} + \OEdgeBN{blue}{blue}{white}{white}\right)\right\rrbracket \geq 0,
    \end{split}
  \end{equation}
  while the equivalent of \eqref{high-degree-edge-cut} is 
  \begin{equation}\label{low-degree-edge-cut}
    \begin{split}
      &\left\llbracket\left(\OEdgeRaR{blue}{blue}{white}{white} - 0.1342\right) \cdot \left(\OEdgeN{blue}{blue}{white} + \OEdgeB{blue}{blue}{white}\right)^2\right. \\
      &+ \left(\frac{1}{2} - \OEdgeR{blue}{blue}{white}\right) \cdot \left(\OEdgeN{blue}{blue}{white} + \OEdgeB{blue}{blue}{white}\right) \cdot \left(\OEdgeRN{blue}{blue}{white}{white} + \OEdgeBR{blue}{blue}{white}{white}\right) \\
      &\left.+ \left(\frac{1}{2} - \OEdgeR{blue}{blue}{white}\right)^2 \cdot \left(\OEdgeNaN{blue}{blue}{white}{white} + \OEdgeBN{blue}{blue}{white}{white}\right)\right\rrbracket
      \geq 0.
    \end{split}
  \end{equation}
  Condition \eqref{high-degree-edges-absent-graph} implies
  \begin{equation}\label{high-degree-edges-absent}
    \Edge{red}{red} \leq 0.02.
  \end{equation}

  Calculations at \url{https://github.com/ignacy123/dense-half-k4-free-flags} show that no graphon can satisfy all of \eqref{degree-blue}, \eqref{degree-red}, \eqref{random-cut},
  \eqref{low-degree-prevalent}, \eqref{high-degree-edges-absent}, \eqref{high-degree-vertex-sides-case-1}, \eqref{high-degree-vertex-cut-case-1},
  \eqref{low-degree-edge-sides}, \eqref{low-degree-edge-cut}, \eqref{degree-cut-sides-case-1} and \eqref{degree-cut-linear-case-1} at once.
  This implies a contradiction similarly to Case~1.

  \textbf{Case 6:}
  In this case, we analyze a $K_4$-free graph $G$ that satisfies
  \begin{equation*}
    \frac{|L(G)|}{n} \geq \frac{1}{2},
  \end{equation*}
  \begin{equation*}
    \frac{e(H(G))}{\binom{n}{2}} \leq 0.02,
  \end{equation*}
  and
  \begin{equation*}
    \begin{split}
      &\left(\frac{\frac{n}{2}}{|L(G)|}\right)^2 e(L(G)) \\
      \geq\ &
      e(H(G)) + \frac{\frac{n}{2} - |H(G)|}{|L(G)|} \cdot e(L(G), H(G)) + \left(\frac{\frac{n}{2} - |H(G)|}{|L(G)|}\right)^2 \cdot e(L(G)).
    \end{split}
  \end{equation*}

  We have already defined all linear inequalities that are necessary to prove a contradiction in this case.

  Calculations at \url{https://github.com/ignacy123/dense-half-k4-free-flags} show that no graphon can satisfy all of \eqref{degree-blue}, \eqref{degree-red}, \eqref{random-cut},
  \eqref{low-degree-prevalent}, \eqref{high-degree-edges-absent}, \eqref{high-degree-vertex-sides-case-1}, \eqref{high-degree-vertex-cut-case-1},
  \eqref{low-degree-edge-sides}, \eqref{low-degree-edge-cut}, \eqref{degree-cut-sides-case-2} and \eqref{degree-cut-linear-case-2} at once.
  This implies a contradiction similarly to Case~1.

  \textbf{Case 7:}
  The exact set of inequalities which $G$ has to satisfy in this case is
  \begin{equation*}
    \frac{|H(G)|}{n} \geq \frac{1}{2},
  \end{equation*}
  \begin{equation*}
    \frac{e(H(G))}{\binom{n}{2}} \leq 0.02,
  \end{equation*}
  and
  \begin{equation*}
    \begin{split}
      &e(L(G)) + \frac{\frac{n}{2} - |L(G)|}{|H(G)|} \cdot e(H(G), L(G)) + \left(\frac{\frac{n}{2} - |L(G)|}{|H(G)|}\right)^2 \cdot e(H(G)) \\
      \geq\ & \left(\frac{\frac{n}{2}}{|H(G)|}\right)^2 \cdot e(H(G)).
    \end{split}
  \end{equation*}

  We have already defined all linear inequalities that are necessary to prove a contradiction in this case.

  Calculations at \url{https://github.com/ignacy123/dense-half-k4-free-flags} show that no graphon can satisfy all of \eqref{degree-blue}, \eqref{degree-red}, \eqref{random-cut},
  \eqref{high-degree-prevalent}, \eqref{high-degree-edges-absent}, \eqref{high-degree-vertex-sides-case-3}, \eqref{high-degree-vertex-cut-case-3},
  \eqref{low-degree-edge-sides}, \eqref{low-degree-edge-cut}, \eqref{degree-cut-sides-case-3} and \eqref{degree-cut-linear-case-3} at once.
  This implies a contradiction similarly to Case~1.

  \textbf{Case 8:}
  Here, we consider the following set of inequalities:
  \begin{equation*}
    \frac{|H(G)|}{n} \geq \frac{1}{2},
  \end{equation*}
  \begin{equation*}
    \frac{e(H(G))}{\binom{n}{2}} \leq 0.02,
  \end{equation*}
  and
  \begin{equation*}
    \begin{split}
      &\left(\frac{\frac{n}{2}}{|H(G)|}\right)^2 e(H(G)) \\
      \geq\ &
      e(L(G)) + \frac{\frac{n}{2} - |L(G)|}{|H(G)|} \cdot e(H(G), L(G)) + \left(\frac{\frac{n}{2} - |L(G)|}{|H(G)|}\right)^2 \cdot e(H(G)).
    \end{split}
  \end{equation*}

  We have already defined all linear inequalities that are necessary to prove a contradiction in this case.

  Calculations at \url{https://github.com/ignacy123/dense-half-k4-free-flags} show that no graphon can satisfy all of \eqref{degree-blue}, \eqref{degree-red}, \eqref{random-cut},
  \eqref{high-degree-prevalent}, \eqref{high-degree-edges-absent}, \eqref{high-degree-vertex-sides-case-3}, \eqref{high-degree-vertex-cut-case-3},
  \eqref{low-degree-edge-sides}, \eqref{low-degree-edge-cut}, \eqref{degree-cut-sides-case-4} and \eqref{degree-cut-linear-case-4} at once.
  This implies a contradiction similarly to Case~1.
\end{proof}

%% file: sections/acknowledgement.tex
The computer-assisted part of this work is based on the flag program originally developed by Bernard Lidický (which was used in \cite{code-bernard}).
The authors are grateful to him for permission to modify the program and to publish the modified version.

%% file: sections/ai.tex
The results in this paper have been obtained without AI assistance.
The paper itself has also been written fully by the authors.
The authors have utilized large language models for linguistic checks and a final correctness check.
Some minor suggestions from AI for the paper itself and related computer-assisted calculations have been implemented.